\documentclass[12pt]{article}
\usepackage{amsmath,amssymb}
\usepackage[]{graphicx}

\def \qed {\hfill \hbox {\rule [-2pt]{3pt}{6pt}}}
\newtheorem{thm}{Theorem}
\newtheorem{theorem}[thm]{Theorem}
\newtheorem{prop}[thm]{Proposition}
\newtheorem{proposition}[thm]{Proposition}
\newtheorem{cor}[thm]{Corollary}

\newtheorem{lem}[thm]{Lemma}

\newtheorem{remark}[thm]{Remark}

\begin{document}
\title{Application of Stochastic Mesh Method to Efficient Approximation of CVA}
 
\author{
Yusuke MORIMOTO
\thanks{
Graduate School of Mathematical Sciences, 
The University of Tokyo, 
Komaba 3-8-1, Meguro-ku, Tokyo 153-8914, Japan, \
Bank of Tokyo Mitsubishi UFJ }
}
\date{}

\maketitle
\begin{abstract}
In this paper, the author considers the numerical computation
of CVA for large systems by Mote Carlo methods.
He introduces two types of stochastic mesh methods
for the computations of CVA. In the first method, stochastic mesh method is used to obtain
the future value of the derivative contracts. In the second method,
stochastic mesh method is used only to judge whether future value of the derivative contracts
is positive or not. He discusses the rate of convergence
to the real CVA value of these methods.
\end{abstract}


Mathematical Subject Classification(2010)   65C05,  60G40


stochastic mesh method, CVA

\section{Introduction}
The credit valuation adjustment (CVA) is, by definition, the difference
between the risk-free portfolio value and the true portfolio value that takes into account
default risk of the counterparty. In other words, CVA is the market value of counterparty credit risk.
After the financial crisis in 2007-2008, it has been widely recognized that even major financial institutions may default.
Therefore, the market participants has become fully aware of counterparty credit risk.
In order to reflect the counterparty credit risk in the price of over-the-counter (OTC) derivative transactions,
CVA is widely used in the financial institutions today.

Although Duffie-Huang \cite{DU} has already introduced the basic idea of CVA  in 1990's, 
several people reconsidered the theory of CVA related to collateralized derivatives (cf. \cite{FT}) and 
also efficient numerical calculation methods appeared(cf. \cite{H}).

There are two approaches to measuring CVA: unilateral and bilateral (cf. \cite{Gre}). 
Under the unilateral approach, it is assumed that the bank that does
the CVA analysis is default-free.
CVA measured in this way is the current market value of future losses due to the counterparty's potential
default. The problem with unilateral CVA is that both the bank and the counterparty require a
premium for the credit risk they are bearing and can never agree on the fair value of the trades in
the portfolio. Therefore, we have to consider not only the market value of the counterparty's default risk, but also the bank's own counterparty credit risk
called debit value adjustment (DVA) in order to calculate the correct fair value. 
Bilateral CVA (it is calculated by netting unilateral CVA and DVA) takes into account the possibility of both the counterparty default and the own default.
It is thus symmetric between the own company and the counterparty, and results in an
objective fair value calculation.

Mathematically, unilateral CVA and DVA are calculated in the same way, and bilateral CVA is the difference of them. So we focus on the calculation of 
unilateral CVA in this paper.

CVA is measured at the counterparty level and there are many assets in the portfolio generally.
Therefore, we have to be involved in the high dimensional numerical problem to obtain the value of CVA.
This is one of the reasons why CVA calculation is difficult.
On the other hand, each payoff usually depends only on a few assets.
We will focus on this property and suggest an efficient calculation methods of CVA 
in the present paper.

Let us consider the portfolio consist of the contracts on one counterparty. 
Let $X^{(m)}(t)$ be ${\bf R}^{N_m}$-valued
stochastic processes $, m=0, 1, \ldots, M.$ We think that $X(t)=(X^{(0)}(t),\ldots, X^{(M)}(t))$ is an underlying process. 
We consider the model that the macro factor is determined by $X^{(0)}(t),$ and the payoff of each derivative at maturity $T_k, k=1,\ldots,K,$ is the form of
$$\sum_{m=1}^M \tilde{F}_{m,k}(X^{(0)}(T_k), X^{(m)}(T_k)).$$
Let $T=T_K$ be the final maturity of all the contracts in the portfolio. Let $\tau$ be the default time of the counterparty, $\lambda(t)$ be its 
hazard rate process, $L(t)$ be the process of loss when the default takes place at time $t$, and $D(t,T)$ be the discount factor process from $t$ to $T$.
We assume that $D(0,t)$ is the function of $X^{(0)}(t)$ and that $L(t), \lambda(t)$ and $\exp(-\int_0^t \lambda(s)ds )$ are the function of $X(t).$   

Let $\tilde{V}_0(t)$ be total value of all contracts in the portfolio at time $t$ under the assumption that counterparty is default free. Then $\tilde{V}_0(t)$ is given by 
$$\tilde{V}_0(t) =E[\sum_{m=1}^M \sum_{k; T_k \geqq t}D(t,T_k)\tilde{F}_{m,k}(X^{(0)}(T_k), X^{(m)}(T_k)) | \mathcal{F}_t],$$
where $E$ denotes the expectation with respect to the risk neutral measure.
Then unilateral CVA on this portfolio is the restructuring cost when the counterparty
defaults.
So unilateral CVA is given by
\begin{align}
\text{CVA}&= E[L(\tau)D(0,\tau)1_{\{\tau <T\}}(\tilde{V}_0(\tau)\vee 0)] \nonumber\\
&=E[\int_0^T L(t)\exp(-\int_0^t \lambda(s)ds ) \lambda(t)D(0,t) (\tilde{V}_0(t) \vee 0)  dt ] \nonumber\\
&=E[\int_0^T L(t)\exp(-\int_0^t \lambda(s)ds ) \lambda(t) (V_0(t) \vee 0)  dt ] \label{defCVA},
\end{align}
where
$$V_0(t) =E[\sum_{m=1}^M \sum_{k; T_k \geqq t}F_{m,k}(X^{(0)}(T_k), X^{(m)}(T_k)) | \mathcal{F}_t],$$
and $F_{m,k}$ is a function such as 
$$F_{m,k}(X^{(0)}(T_k), X^{(m)}(T_k))=D(0,T_k)\tilde{F}_{m,k}(X^{(0)}(T_k), X^{(m)}(T_k)).$$
Since $L(t)\exp(-\int_0^t \lambda(s)ds ) \lambda(t)$ is a function of $X(t),$ we denote it by $g(t,X(t)).$ Then CVA is given by the following form. 
\begin{align}
\text{CVA}&=E[\int_0^T g(t,X(t))  (E[\sum_{m=1}^M \sum_{k; T_k \geqq t}F_{m,k}(X^{(0)}(T_k), X^{(m)}(T_k)) | \mathcal{F}_t] \vee 0)  dt ] \label{defCVA2}.
\end{align}

Now we prepare the mathematical setting. 
Let $M \geqq 1$ be fixed, ${N}_{m} \geqq 1, m=1, \cdots, M$, 
$N= N_0+ \cdots + N_M,$ $ \tilde{N}_m=N_0+N_m,$ and $\tilde{N}=\max_{m=1\ldots,M}\tilde{N}_m$.\\
Let $W_0$ $ = \{ w\in C([0,\infty );{\bf R}^d); \; w(0) = 0 \} ,$
${\cal F}$ be the Borel algebra over $W_0$ 
and $\mu$ be the Wiener measure on $(W_0,{\cal F}).$
Let $B^i:[0,\infty )\times W_0 \to {\bf R},$ $i=1,\ldots ,d,$ be given 
by $B^i(t,w) =w^i(t),$ $(t,w)\in [0,\infty )\times W_0.$
Then $ \{ (B^1(t), \dots ,B^d(t) ; t \in [0,\infty ) \}$ 
is a $d$-dimensional Brownian motion.
Let $B^0(t) = t,$ $t \in [0,\infty ).$

\noindent Let $V^{(0)}_i$ $ \in C^{\infty}_b({\bf R}^{{N}_0}; {\bf R}^{{N}_0}),$ 
$V^{(m)}_i \in C^{\infty}_b({\bf R}^{{N}_0}\times {\bf R}^{{N}_m};{\bf R}^{{N}_m}), i =0, \cdots, d, m=1, \cdots, M.$
Here $C^{\infty}_b({\bf R}^m;{\bf R}^n)$ denotes 
the space of ${\bf R}^n$-valued smooth functions defined 
in ${\bf R}^m$ whose derivatives of any order are bounded.
We regard elements in $C^{\infty}_b({\bf R}^n;{\bf R}^n)$ 
as vector fields on ${\bf R}^n.$

\noindent Now let us consider the following Stratonovich stochastic differential equations.
\begin{equation}
{X}^{(0)}(t,x_0) = {x}_0 + \sum_{i=0}^{d} \int_{0}^{t} V^{(0)}_{ i}({X}^{(0)}(s,{x}_0))\circ dB_i(s),
\label{eq:SDE0}
\end{equation}
\begin{align}
{X}^{(m)}(t,\tilde{x}_m) = {x}_m +  \sum_{i=1}^{d} \int_{0}^{t} V^{(m)}_{i}({X}^{(0)}(s,{x}_0), {X}^{(m)}(s,\tilde{x}_m))\circ dB_i(s),
\label{eq:SDEk}
\end{align}
where $x_m \in {\bf R}^{N_m}, \tilde{x}_m=(x_0, x_m) \in {\bf R}^{{N}_0}\times {\bf R}^{{N}_m}, m=1,\ldots,M.$\\
 Let $\tilde{X}^{(m)}(t,\tilde{x}_m) = ({X}^{(0)}(t,{x}_0), {X}^{(m)}(t,\tilde{x}_m))$ and $\tilde{V}^{(m)}_ i \in C^{\infty}_b({\bf R}^{{N}_0}\times {\bf R}^{{N}_k};{\bf R}^{{N}_0}\times{\bf R}^{{N}_k}), i =0, \cdots, d, \ m = 1,\cdots, M$ be \begin{align*}
 \tilde{V}^{(m)}_{i}(\tilde{x}_m)= 
 \begin{pmatrix}
 V^{(0)}_{ i}(x_0)\\
 V^{(m)}_{ i}(\tilde{x}_m)
 \end{pmatrix} .
 \end{align*}
Then we have
\begin{align}
\tilde{X}^{(m)}(t,\tilde{x}_m) = \tilde{x}_m + \sum_{i=0}^{d} \int_{0}^{t} \tilde{V}^{(m)}_i(\tilde{X}^{(m)}(t,\tilde{x}_m) )\circ d B_i(s).
\end{align}
There is a unique solution $\tilde{X}^{(m)}(t,\tilde{x}_m) $ 
to this equation. 
Then $X(t,x), x\in {\bf R}^N$ also satisfies the solution to the following Stratonovich stochastic differential equation. 
\begin{equation}
X(t,x) = x + \sum_{i=0}^{d} \int_{0}^{t} {V}_i(X(s,x))\circ d {B}^i(s),
\label{eq:SDE}
\end{equation}
where ${V}_i, i =1,\ldots, d$ is
\begin{align*}
 {V}_{i}(x)= 
 \begin{pmatrix}
 V^{(0)}_{ i}(x_0)\\
  V^{(1)}_{ i}(\tilde{x}_1)\\
 \vdots \\
 V^{(M)}_{ i}(\tilde{x}_M)
 \end{pmatrix} .
 \end{align*}

We assume that vector fields ${V}_i, i=1,\ldots,d,$ satisfy condition (UFG) stated in the section \ref{vec}. 
Let $E_m $ be defined by (\ref{defE}) in Section \ref{vec}. By \cite{KM}, if $\tilde{x}_m \in E_m,$ 
the distribution law of $\tilde{X}^{(m)}(t,\tilde{x}_m)$ under $\mu$ has a smooth density function
$p^{(m)}(t,\tilde{x_m},\cdot ):{\bf R}^{\tilde{N}_m}\to [0,\infty )$ for $t>0, m=1, \ldots, M.$\\
Let $x^{*} =(x_0^*, \ldots, x_M^*) \in {\bf R}^N.$ We assume that the underlying asset process is $X(t)=X(t, x^*)$.
We also assume that
$$\tilde{x}^*_m=(x_0^*, x_m^*) \in E_m, m=1,\ldots, M.$$
Let $\hat{{\cal D}}({\bf R}^n)$ denotes the space of functions on ${\bf R}^n$ given by
 $$\hat{{\cal D}}({\bf R}^n) =\{ f \in C^2({\bf R}^n); \|\frac{\partial^{\alpha}f}{\partial x^{\alpha}}\|_{\infty} <\infty, \text{for  } 1 \leqq |\alpha| \leqq 2 \},$$
 where $\|f\|_{\infty}=\sup \{|f(x)|; x\in {\bf R}^n\}$.\\
${Lip}({\bf R}^n)$ denotes the space of Lipschitz continuous functions on ${\bf R}^n,$  and
we define a semi-norm $\|\cdot \|_{Lip}$ on $Lip({\bf R}^n)$ by
$$\|f\|_{{Lip}}=\sup_{x,y\in {\bf R}^n, x \neq y} \frac{|f(x)-f(y)|}{|x-y|},
\quad f\in {Lip}({\bf R}^n).$$
Let ${\cal M}({\bf R}^n)$ be the linear subspace of ${Lip}({\bf R}^n)$  spanned by
$\{f \vee g;  f, g \in {\cal D}({\bf R}^n)\}.$\\
\\
We define linear operators 
$P_{t}:{Lip}({\bf R}^N) \to {Lip}({\bf R}^N), t \geqq 0,$ by
$$
(P_{t}f)(x) = E^{\mu}[f(X(t,x))], \quad f\in {Lip}({\bf R}^N),
$$
and $P^{(m)}_{t}:{Lip}({\bf R}^{\tilde{N}_m}) \to {Lip}({\bf R}^{\tilde{N}_m}), t \geqq 0 , m=1, \ldots,M,$ by
$$
(P^{(m)}_{t}f)(\tilde{x}_m) = E^{\mu}[f(\tilde{X}^{(m)}(t,\tilde{x}_m))], \quad f \in {Lip}({\bf R}^{\tilde{N}_m}).
$$
We remind that $L(t)\exp(-\int_0^t\lambda(s)ds)\lambda(t) $ is represented by
$$
L(t)\exp(-\int_0^t\lambda(s)ds)\lambda(t) = g(t, X(t,x^*)).
$$
We assume that $g: [0,T] \times {\bf R}^N \to [0,\infty)$ satisfies the following two conditions.\\
(1) $g(t,x)$ is differentiable in $t$ and there is an integer $n_1,$ and a constant $C_1 >0$ such that
$$\sup_{t\in[0,T]}| \frac{\partial}{\partial t}g(t,x)|  \leqq C_1(1+|x|^{n_1})
,\quad x \in {\bf R}^N.$$
(2) $g(t,x)$ is $2$-times coninuously differentiable in $x$ and there is an integer $n_2,$ and a constant $C_2 >0$ such that
$$
\sup_{t\in[0,T]} |\frac{\partial^{\alpha}}{\partial x^{\alpha}}g(t,x)| \leqq C_2(1+|x|^{n_2}),\quad x \in {\bf R}^N
$$
for any multi index $|\alpha|\leqq 2.$
\\
We assume that a discounted payoff functions $F_{m,k}, m=1,\ldots, M, k=1,\ldots,K$ in equation (\ref{defCVA2}) 
belong to ${\cal M}({\bf R}^{\tilde{N}_m}).$
Under the assumptions above, CVA $c_0$ is given by 
\begin{align}
c_0 = E^{\mu}[\int_{0}^T \{ g(t,X(t,x^*)) 
E^{\mu}[ \sum_{m=1}^M \sum_{k:T_k\geqq t} F_{m,k}(\tilde{X}^{(m)}(T_k, \tilde{x}^*_m)) |\mathcal{F}_t] \vee0\}dt ].
\end{align}
We will introduce numerical calculation methods by Monte Carlo simulation for $c_0$.\\
\\
Let $(\Omega ,{\cal F},P)$ be a probability space,
and $X_{\ell }:[0,\infty )\times \Omega \to {\bf R}^N,$ 
$\ell =1,2,\ldots ,$ be continuous stochastic processes such that 
each probability law on $C([0,\infty );{\bf R}^N)$ of
$X_{\ell}(\cdot )$ under $P$ is the same as  that of $X(\cdot ,{x}^*)$  for all $\ell$ $=1,2,\ldots ,$ and that
$\sigma \{ X_{\ell}(t); \; t\geqq 0\},$ $\ell =1,2,\ldots ,$ are independent.\\
Let us define projections $\pi_m: {\bf R}^N \to {\bf R}^{\tilde{N}_m}, m=1,\ldots,M,$ by
$ \pi_m(x)= \tilde{x}_m=(x_0,x_m),$ 
and define $\varepsilon_0 >0$ by
$\varepsilon_0 =\min_{1\leqq k \leqq K}(T_k-T_{k-1}).$
We define random linear operators (stochastic mesh operators)
 $Q_{t,T_k,\varepsilon}^{(m)}=Q_{t,T_k,\varepsilon}^{(m,L,\omega)}, 0 \leqq t \leqq T, 0 < \varepsilon <\varepsilon_0,$ 
on $Lip({\bf R}^{\tilde{N}_m})$ by
\begin{align*}
&(Q_{t,T_k,\varepsilon}^{(m,L,\omega)}f)(\tilde{x}_m)
=
\begin{cases}
\frac{1}{L} \sum_{\ell=1}^L \frac{f(X_{\ell}^{(m)}(T_k))p^{(m)}(T_k-t, \tilde{x}_m, \pi_m( X_{\ell}(T_k)))}{q_{t,T_k}^{(m,L,\omega)}(\pi_m (X_{\ell}(T_k)))}, 
\quad 0 \leqq t < T_k-\varepsilon ,\\
f(\tilde{x}_m), \quad T_k - \varepsilon \leqq t \leqq T_k, \\
0, \quad T_k < t \leqq T .
\end{cases}\\
&\text{where} \quad q_{t,T_k}^{(m,L,\omega)}(\tilde{y}_m)=\frac{1}{L}\sum_{\ell=1}^L p^{(m)}(T_k-t, \pi_m(X_{\ell}(t)), \tilde{y}_m)).
\end{align*}

Let $\Pi$ denotes the set of partitions $\Delta=\{t_0,t_1,\ldots, t_n\}$ such that 
$0=t_0 < t_1 < \ldots < t_n= T$ and that   
$\{T_k; k=1,\ldots,K\} \subset \Delta.$
Let $|\Delta|=\max_{i=1,\ldots,n}{(t_{i+1}-t_i)}.$
We define estimators $\hat{c}_i=\hat{c}_i(\varepsilon,\Delta, L) : \Omega \to {\bf R}, i=1,2,$ in the following.
\begin{align}
&\hat{c}_1(\varepsilon, \Delta, L)(\omega) \nonumber\\
=&\frac{1}{L} \sum_{\ell=1}^L \sum_{i=0}^{n-1} (t_{i+1}-t_i) g(t_i,X_{\ell}(t_i))
(\sum_{m=1}^M \sum_{k:T_k\geqq t_{i+1}} (Q_{t_i,T_k,\varepsilon}^{(m,L,\omega)}F_{m,k})(\pi_k(X_{\ell}(t_i)))\vee0), \label{app1}\\
\text{and}\nonumber \\
&\hat{c}_2(\varepsilon, \Delta, L)(\omega)\nonumber\\
=&\sum_{i=0}^{n-1}(t_{i+1}-t_{i}) E^{\mu}[g(t_i,X(t_i, x^*)) ( \sum_{m=1}^M \sum_{k:T_k\geqq t_{i+1}} F_{m,k}(\pi_k X(T_k, x^*))) \nonumber \\
& \qquad \qquad \qquad \qquad \qquad \qquad \qquad \qquad \times 1_{\{\sum_{m=1}^M \sum_{k:T_k\geqq t_{i+1}} (Q_{t_i,T_k,\varepsilon}^{(m,L,\omega)}F_{m,k})(\pi_k X(t_i,x^*))) \geqq 0\}}]. \label{app2}
\end{align}
Our main results are following.
\begin{theorem} \label{main2}
Let $\alpha_0 = (1+\delta)(\tilde{N}+1)\ell_0/4 \vee 1.$ Let $\{\varepsilon_L\}_{L = 1}^{\infty} \subset (0, \varepsilon_0)$ be a sequence and suppose
 that there is a constant $C_0 \in (0, \infty)$ such that $ \varepsilon_L \leqq C_0 L^{-\frac{1+ \delta}{2(1+\alpha_0)}}, L\geqq 1.$
 Then there exists a constant $C_1 \in (0,\infty)$ such that
$$
E^P[|\hat{c}_1(\varepsilon_L, \Delta, L)-c_{0}|] \leqq C_1 (L^{-\frac{1}{1+\alpha_0}}+|\Delta|)
$$
for any $L \geqq 1$ and $\Delta \in \Pi.$
\end{theorem}

\begin{theorem}\label{main1} 
Let $\alpha_1 = (1+\delta)(\tilde{N}+1)\ell_0/2 \vee 1,$ and
let $\{\varepsilon_L\}_{L = 1}^{\infty} \subset (0, \varepsilon_0)$ be a sequence such that there is a constant $C_0 \in (0,\infty),$ 
such that $\varepsilon_L \leqq C_0 L^{-\frac{1+ \delta}{2\alpha_1+1}}, L \geqq 1.$ 
Suppose that there are constants $ \gamma \in (0,1]$  and $C_{\gamma} \in (0, \infty)$ such that
 $$\sup_{\Delta} \sum_{i=0}^{n-1} (t_{i+1}-t_i)){\mu} (| \sum_{m=1}^M\sum_{k:T_k\geqq t_{i+1}} (P_{T_k-t_i}^{(m)}F_{m,k})(\pi_m X(t_i,x^*)))| \leqq \theta ) \leqq C_{\gamma}\theta^{\gamma}$$
for all $ \theta \in (0,1].$\\
Then then there exists a constant $C_1\ \in (0,\infty)$ and $\tilde{\Omega}(L) \in \mathcal{F},$  $L\geqq 1 ,$ such that
$$P(\tilde{\Omega}(L)) \to 1, \ L\to \infty,$$
and
$$
1_{\tilde{\Omega}(L)}|\hat{c}_2(\varepsilon_L, \Delta, L)-c_{0}|
\leqq  C_1
(L^{-(\frac{1}{2}+ \frac{(1-\delta)}{2 \alpha_1 +1}) \frac{1+\gamma}{2+\gamma}}+|\Delta|)
$$
for any $L\geqq 1,$ and $\Delta \in \Pi.$
\end{theorem}

\begin{remark}\label{actual}
Let $\tilde{\Omega}'(L)$ be
$$ \tilde{\Omega}'(L)=\left \{\omega \in \Omega; |\hat{c}_1(\varepsilon_L,\Delta,L) - c_0 | \leqq CL^{-\frac{1-\delta}{1+\alpha_0}} \right\}.$$
Then by Theorem \ref{main2}, we see that
$$
P(\tilde{\Omega}'(L)) \to 1,  L\to \infty,
$$
and
$$1_{\tilde{\Omega}(L) }|\hat{c}_1(\varepsilon_L,\Delta,L)-c_0| \leqq CL^{-\frac{1-\delta}{1+\alpha_0}}.
$$
Theorem \ref{main1} shows that the estimation of $\hat{c}_2$ may be better than $\hat{c}_1$. 
\end{remark}

We can compute the estimators $\hat{c}_i, i=1,2,$ practically in the following way.\\  
\noindent First, we generate a family of independent paths 
$${\bf X}_1 = \{X_{\ell }(t); 0 \leqq t \leqq T, \ell =1,2,\ldots ,L\}.$$
Next, by using ${\bf X}_1$, we compute
$$(Q_{t_i,T_k,\varepsilon}^{(m,L,\omega)}F_{m,k})(\pi_k(X_m(t_i))), \text{for every}\  k \  \text{such that} \ T_k > t_i .$$
Then our estimator $\tilde{c}_1$ is
$$\tilde{c}_1=\frac{1}{L} \sum_{\ell=1}^L \sum_{i=0}^{n-1} g(t_i, X_{\ell}(t_i))( \sum_{k:T_k\geqq t_{i+1}} (Q_{t_i,T_k,\varepsilon}^{(m,L,\omega)}F_{m,k})(\pi_k(X_{\ell}(t_i)))\vee0)(t_{i+1}-t_i). $$
We used the same paths for Monte Carlo simulation and construction  of Stochastic mesh operator. 

For $\tilde{c}_2$, we generate another independent family of independent paths $${\bf X}_2=\{X^{'}_{m }(t); 0 \leqq t \leqq T,  m =1,2,\ldots ,M\},$$
and we compute 
\begin{multline*}
\tilde{c}_2= \frac{1}{M} \sum_{m=1}^M \sum_{i=0}^{n-1} g(t_i, X'_m(t_i)) \{ \sum_{k:T_k\geqq t_{i+1}}F_{m,k}(X_m^{'(k)}(T_k))\}\\
 \times 1_{\{\sum_{k:T_k\geqq t_{i+1}} (Q_{t,T,\varepsilon}^{(m,L,\omega)}F_{m,k})(\pi_k(X^{'}_m(t_i)))) \geqq 0\}}(t_{i+1}-t_i).
 \end{multline*}

In the above computations of $\tilde{c}_1$ and $\tilde{c}_2$, we use the values of $X_{\ell}(t_i), t_i \in \Delta$, only.\\
As for the computation of $\tilde{c}_2,$  we do not use $Q_{t,T,\varepsilon}^{(m,L,\omega)}F_{m,k}$ explicitly.
We use $Q_{t,T,\varepsilon}^{(m,L,\omega)}F_{m,k}$ only to judge whether $(P^{(k)}_{T_k-t} F_{m,k})(\pi_k(X^{'}_m(t_i))))>0$ or not. 
 So the approximation has no error when the signs of  $(Q_{t,T,\varepsilon}^{(m,L,\omega)} F_{m,k})(\pi_k{X}(t)))$ and $(P^{(k)}_{T_k-t} F_{m,k})(\pi_k(X^{'}_m(t_i)))$ are the same, even if there are large differences between them.

\section{Structure of Vector Fields}\label{vec}
Let ${\cal A}=\bigcup_{k = 1}^{\infty} \{ 0,1,\ldots , d \}^k$
and ${\cal A}^{*}= {\cal A} \setminus \{ 0\} .$ 
For $\alpha \in {\cal A}$, 
Let $|\alpha | = k$ if 
$\alpha $ $= (\alpha^1,\ldots ,\alpha^k)$
$\in \{ 0,1,\ldots , d \}^k ,$ and 
let $ \parallel \alpha \parallel $ 
$ = |\alpha | + \mbox{card}\{ 1 \leqq i \leqq |\alpha | ; 
\; \alpha^i = 0 \}.$ 
Also, for each $m\geqq 1,$ ${\cal A}^{*}_{\leqq m}  
=\{ \alpha \in {\cal A}^{*}; \; \parallel \alpha \parallel \leqq m \} .$

We define vector fields $V_{[\alpha ]},$ $\alpha \in {\cal A},$ 
inductively by
$$
\qquad V_{[i]} = V_i, 
\quad i = 0,1,\ldots ,d,
$$
$$
V_{[\alpha *i]} = [V_{[\alpha ]},V_i], 
\qquad i = 0,1,\ldots ,d.
$$
Here $\alpha *i$ 
$= (\alpha^1,\ldots ,\alpha^k,i)$
for $\alpha $ $= (\alpha^1,\ldots ,\alpha^k)$ and $i = 0,1,\ldots ,d.$

We assume that a system of vector fields $\{ V_i ; i=0,1,\ldots ,d \}$ 
satisfies the following condition (UFG).

\noindent
(UFG) There is an integer $\ell_0$ 
and there are functions $\varphi_{\alpha,\beta} \in C^{\infty}_b({\bf R}^N),$ 
$\alpha \in {\cal A}^{*},$ $\beta \in {\cal A}^{*}_{\leqq \ell_0},$ 
satisfying the following.
$$
{\displaystyle V_{[\alpha ]} 
= \sum_{\beta \in {\cal A}^{*}_{\leqq \ell_0}}
\varphi_{\alpha,\beta}V_{[\beta ]}},
\qquad 
\alpha \in {\cal A}^{*}.
$$

\begin{prop}
A system of vector fields $\{ \tilde{V}^{(m)}_{i}; i=0,1,\ldots , {d} \}$  
also satisfies the $(UFG)$ condition.
\end{prop}
 {\it Proof.}
We prove following by induction on $|\alpha|.$
\begin{align} \label{inductive}
{V}_{[\alpha ]} (f \circ \pi_m) = (\tilde{V}^{(m)}_{[\alpha ]}f)\circ \pi_m, \quad f \in C_b^{\infty}({\bf R}^{\tilde{N}_m}),
 \end{align}
 for any $\alpha\in \mathcal{A}$ and $m=1,\ldots,M.$ \\
It is trivial in the case of $|\alpha|=1.$
 By the assumption for induction,
$$
{V}_{[\alpha * i]} (f \circ \pi_m) =({V}_{[\alpha]} V_i- V_i {V}_{[\alpha]}) (f \circ \pi_m)
$$
$$
={V}_{[\alpha]} ((\tilde{V}^{(m)}_{i} f)\circ \pi_m )- V_i ((\tilde{V}^{(m)}_{[\alpha]}f )  \circ \pi_m)
$$
$$
=(\tilde{V}^{(m)}_{[\alpha]}(\tilde{V}^{(m)}_{i} f)) \circ \pi_m - (\tilde{V}^{(m)}_{i} (\tilde{V}^{(m)}_{[\alpha]}f ))  \circ \pi_m.
$$
So we have (\ref{inductive}).
From (UFG) condition, we have
$$ {V}_{[\alpha ]} (f \circ \pi_m)= 
\sum_{ \beta \in {\cal A}^{*}_{\leqq \ell_0}} \varphi_{ \alpha, \beta } 
{V}_{[\beta]} (f \circ \pi_m)$$
$$=\sum_{ \beta \in {\cal A}^{*}_{\leqq \ell_0}} \varphi_{ \alpha, \beta } 
(\tilde{V}_{[\beta]}^{(m)}f) \circ \pi_m.$$
 Let $j_m: {\bf R}^{\tilde{N}_m} \to {\bf R}^{N} $ be
 $$j_m(\tilde{x}_m) = (x_0,0\ldots , 0,x_m,0,\ldots,0).$$
Then
$$
\tilde{V}^{(m)}_{[\alpha]}f=(V_{[\alpha]}^{(m)}f)\circ \pi_m \circ j_m
$$
$$
=(\sum_{ \beta \in {\cal A}^{*}_{\leqq \ell_0}} \varphi_{ \alpha, \beta } 
(\tilde{V}_{[\beta]}^{(m)}f) \circ \pi_m) \circ j_m
=\sum_{ \beta \in {\cal A}^{*}_{\leqq \ell_0}} (\varphi_{ \alpha, \beta }\circ j_m)
\tilde{V}_{[\beta]}^{(m)}f.
$$
So we have our assertion.
\qed \\
\\
Let $A_m(\tilde{x}_m)= (A_m^{ij}(\tilde{x}_m))_{i,j=1,\ldots ,\tilde{N}_m},$ $t>0,$ $\tilde{x}_m\in {\bf R}^{\tilde{N}_m},$ 
be a $\tilde{N}_m \times \tilde{N}_m$ symmetric matrix given by
$$
A_m^{ij}(\tilde{x}_m) 
= \sum_{\alpha \in {\cal A}^{*}_{m, \leqq \ell_0}}
\tilde{V}^i_{m,[\alpha ]}(\tilde{x}_m)\tilde{V}^j_{m,[\alpha ] }(\tilde{x}_m),
\qquad i,j=1,\ldots , \tilde{N}_m.
$$
Let $h_m(\tilde{x}_m) = \det A_m(\tilde{x}_m), \tilde{x}_m \in {\bf R}^{ \tilde{N}_m}$ 
and 
\begin{align}
E_m= \{ \tilde{x}_m \in {\bf R}^{\tilde{N}_m};\; h_m(\tilde{x}_m)>0\} . 
\label{defE}
\end{align}
By \cite{KS2}, we see that if $\tilde{x}_m\in E_m,$ 
the distribution law of $\tilde{X}^{(m)}(t,\tilde{x}_m)$ under $\mu$ has a smooth density function
$p^{(m)}(t,\tilde{x}_m,\cdot ):{\bf R}^{\tilde{N}_m}\to [0,\infty )$ for $t>0.$
Moreover, by \cite{KM} we see  that $\int_{E_m} p^{(m)}(t,\tilde{x}_m,\tilde{y}_m)dy =1,$ $\tilde{x}_m\in E_m.$ 
We have $p_m(t,\tilde{x}_m,y)=0, y\in E_m^c$ by \cite{KM}.\\

\section{Prepareations}
In this section, we use the notation in \cite{K}.
We have the following Lemma similarly to the proof of \cite{K} Lemma 8 (3).  
\begin{lem} \label{rev}
For any $\Phi \in {\bf D}_{\infty-}^1, \alpha \in \mathcal{A}^*_{\leqq \ell_0},$ let
$$(D^{(\beta)} \Phi)(t,x)=(D \Phi(t,x), k^{\beta}(t,x))_H$$
and
\begin{align*}
\Phi_{\alpha}(t,x)=&\sum_{\beta \in \mathcal{A}^*_{\leqq \ell_0}} t^{-\|\alpha\| / 2}\{ -D^{(\beta)} \Phi(t,x)M_{\alpha \beta}^{-1}(t,x) \\
&-(\sum_{\gamma_1 \gamma_2 \in \mathcal{A}^*_{\leqq \ell_0}} \Phi(t,x) M_{\alpha \gamma_1}^{-1}(t,x)) D^{(\beta )}
M^{\gamma_1 \gamma_2}(t,x) M_{\gamma_1 \beta}^{-1}(t,x) \\
&+ \Phi(t,x) M_{\alpha \beta}^{-1}(t,x)) D^* k^{\beta}(t,x) \}, \quad t>0, x \in {\bf R}^{N}.
\end{align*}
Then 
$$
E^{\mu}[\Phi(t,x) (V_{[\alpha]}f )(X(t,x))] = t^{-\|\alpha\| / 2}E^{\mu}[\Phi_{\alpha}(t,x) f(X(t,x))],
$$
and
$$\sup_{t\in[0,T], x \in {\bf R}^N, p\in (1, \infty)} E[|\Phi_{\alpha}(t,x)|^p] < \infty.$$
\\
\end{lem}

Let $\varphi$ be a smooth function such that
\begin{align}
&\varphi(z)=
\begin{cases}
1 , z\geqq 1\\
0, z<0,
\end{cases}\\
&\varphi' (z) \geqq 0.
\end{align} 
Let $\varphi_m(z)=\varphi(mz)$ and $\bar{\varphi}$ be 
\begin{align} \label{defPhi}
\bar{\varphi}_m(z)=\int_{0}^z \varphi_m(z')dz'.
\end{align}
Then for any $z \in {\bf R},$
$$\bar{\varphi}_m(z) \to z \vee 0 , \quad m\to \infty.$$ 

\begin{lem}\label{rev2}
If $\Phi \in {\bf D}_{\infty -}^1,$ then $|\Phi| \in {\bf D}_{\infty -}^1.$
\end{lem}
{\it Proof.}
Let $\bar{\psi}_m(z) = \bar{\varphi}_m(z)+\bar{\varphi}_m(-z).$ Then for any $z \in {\bf R},$
$$\bar{\psi}_m(z) \to |z|  , \quad m\to \infty,$$
and $|\bar{\psi}'_m(z)| \leqq 1.$
We have
\begin{align*}
&D(\bar{\psi}_m(\Phi(t,x))) = \bar{\psi}'_m(\Phi(t,x))D\Phi(t,x)\\
=&(\varphi_m(\Phi(t,x))-\varphi_m(-\Phi(t,x)))D\Phi(t,x), \quad m\geqq 1.
\end{align*}
Then $\{D(\bar{\psi}_m(\Phi(t,x)))\}_{m=1}^{\infty}$ is a Cauchy sequence in $L^p(W_0,\mathcal{L}_{(2)}^1(H;{\bf R})), p>1,$
because 
\begin{align*}
\|D(\bar{\psi}_m(\Phi(t,x))) - D(\bar{\psi}_n(\Phi(t,x))) \|_H \leqq 1_{\{|\Phi(t,x)|\in [0,1/m] \}}\|D\Phi(t,x)\|_H, \quad n\geqq m\geqq 1.
\end{align*}
Because $D: {\bf D}_p^1\to L^p(W_0,\mathcal{L}_{(2)}^1(H;{\bf R}))$ is a closed operator,
we have $|\Phi(t,x)|\in {\bf D}_p^1,$ for any $p>1$. So we have the assertion.
\qed
\\

Let us denote $\|\nabla F\|_{\infty}=\sup_{x \in {\bf R}^N}|(\frac{\partial F}{\partial x_1}(x), \ldots, \frac{\partial F}{\partial x_N}(x))|, \quad F \in C^{\infty}({\bf R}^N).$
\begin{lem}\label{discLem2}
Let $T>0$. Then there exists a $C>0$ such that 
\begin{align*}
&E[|g(t,X(t,x^*))(P_{T-t}F)(X(t,x^*)) \vee 0 - g(s,X(s,x^*)(P_{T-s}F)(X(s,x^*)) \vee 0|]\\
&\leqq C \|\nabla F\|_{\infty} \int_{s}^t(r^{-1/2}+(T-r)^{-1/2})dr ,
\end{align*}
for any $ F \in C_b^{\infty}({\bf R}^{N})$ and any $0< s < t< T$. 
\end{lem}
{\it Proof.}
Let $\{M(t)\}_{0\leqq t \leqq T}$ be
$$M(t) = E^{\mu}[F(X(T,x^*))|\mathcal{F}_t]=(P_{T-t}F)(X(t,x^*)).$$
$\{M(t)\}_{0\leqq t \leqq T}$ is a $\{\mathcal{F}_t\}_{t\geqq 0}$-martingale.
Let $Y(t)=g(t,X(t,x^*)), 0 \leqq t \leqq T.$

Let 
$$L_t = \frac{\partial}{\partial t} + V_0+\frac{1}{2}\sum_{i=1}^d V_i^2.$$
By It\^o formula, 
\begin{align*}
Y(t)\bar{\varphi}_m(M(t))&=Y(s)\bar{\varphi}_m(M(s))+\int_s^t Y(r)\varphi_m(M(r))dM(r)\\
&+\frac{1}{2}\int_s^t Y(r)\varphi_m'(M(r))d\langle M \rangle(r)\\
&+\int_s^t \bar{\varphi}_m(M(r))dY(r)+\int_s^t d\langle Y, \bar{\varphi}_m(M) \rangle(r).
\end{align*}
Note that,
$$M(t)=M(s) + \sum_{j=1}^d \int_s^t V_j (P_{T-r}F)(X(r,x^*))dB^j(r),$$
$$\langle M \rangle(t)=\langle M \rangle(s) + \sum_{j=1}^d \int_s^t(V_j (P_{T-r}F)(X(r,x^*)))^2dr,$$
$$Y(t)=Y(s) + \sum_{j=1}^d \int_s^t  (V_jg)(X(r,x^*)) dB^j(r) +\int_s^t (L_t g)(X(r,x^*))dr,$$
and
$$\langle Y, \bar{\varphi}_m(M) \rangle(t)=\langle Y, \bar{\varphi}_m(M) \rangle(s) + \sum_{j=1}^d \int_s^t (V_jg)(X(r,x^*))(V_j (P_{T-r}F))(X(r,x^*))dr.$$
So we have
$$E^{\mu}[|Y(t)\bar{\varphi}_m(M(t)) - Y(s)\bar{\varphi}_m(M(s))|]$$
$$=\frac{1}{2}\sum_{j=1}^d \int_s^t E^{\mu}[|Y(r)\varphi_m'( (P_{T-r}F)(X(r,x^*)) ) \left(V_j (P_{T-r}F)(X(r,x^*))\right)^2|] dr$$
$$+\int_s^t E^{\mu}[|\bar{\varphi}_m(M(r)) (L_t g)(X(r,x^*))|]dr+ \sum_{j=1}^d \int_s^t E^{\mu}[|(V_jg)(V_j (P_{T-r}F))(X(r,x^*))|]dr.$$
Now by the definition of $\varphi_m$ and $g$,  we have,
$$\int_s^t E^{\mu}[|\bar{\varphi}_m(M(r)) (L_t g)(X(r,x^*))|]dr \leqq \int_s^t  E^{\mu}[|M(r)|^2]^{1/2}E^{\mu}[|(L_t g)(X(r,x^*))|^2]^{1/2}dr$$
$$ \leqq \sup_{t\in[0,T]}E^{\mu}[|(L_t g)(X(r,x^*))|^2]^{1/2} \int_s^t  E^{\mu}[|M(r)|^2]^{1/2}dr.$$
By Burkholder's inequality,
$$\int_s^t  E^{\mu}[|M(r)|^2]^{1/2}dr \leqq E^{\mu}[\langle M \rangle_t]^{1/2} (t-s)
 \leqq \sup_{r\in (s,t)} \|V_j (P_{T-r}F)\|_{\infty}(t-s).$$
 On the other hand, we have,
$$\sum_{j=1}^d \int_s^t E^{\mu}[|(V_jg)(X(r,x^*))(V_j (P_{T-r}F))(X(r,x^*))|]dr$$
$$ \leqq \|V_j (P_{T-r}F)\|_{\infty} \sum_{j=1}^d \int_s^t E^{\mu}[|(V_jg)(X(r,x^*))|]dr$$
$$\leqq \sum_{j=1}^d \sup_{r\in (s,t)} \|V_j (P_{T-r}F)\|_{\infty} (t-s),$$
for any $ F \in C_b^{\infty}({\bf R}^{N})$ and any $0< s < t< T$.\\ 
On the other hand
$$ \varphi_m' \left( (P_{T-r}F)(x^*) \right) \left(V_j (P_{T-r}F)(x^*)\right)^2$$
$$=\left(V_j (\varphi_m \circ (P_{T-r}F))\right) (x^*) \left(V_j (P_{T-r}F)\right)(x^*) $$
$$= V_j \left( \varphi_m \circ ( P_{T-r}F)(x^*)  V_j (P_{T-r}F)(x^*) \right) 
-  \varphi_m \circ (P_{T-r}F)(x^*)  V_j^2(P_{T-r}F) (x^*) . $$
Notice that  $\varphi_m' \geqq 0,$ we have
$$ E^{\mu}[|Y(r)\varphi_m' \left( (P_{T-r}F)(X(r,x^*)) \right) \left(V_j (P_{T-r}F)(X(r,x^*))\right)^2|]$$
$$= E^{\mu}[|Y(r)| \varphi_m' \left( (P_{T-r}F)(X(r,x^*)) \right) \left(V_j (P_{T-r}F)(X(r,x^*))\right)^2]$$
$$=  I_{1,j}(r,f)-I_{2,j}(r,f),$$
where
$$I_{1,j}( r,F) =E^{\mu} [ |g(r, X(r,x^*))| V_j \left( \varphi_m \circ ( P_{T-r}F)  V_j (P_{T-r}F) \right) (X(r,x^*))], $$
$$I_{2,j}( r,F) = E^{\mu} [  |g(r, X(r,x^*))|  \varphi_m \circ (P_{T-r}F)(X(r,x^*)) V_j^2(P_{T-r}F) (X(r,x^*))]. $$
Let $\Phi_{g}(r,x)=|g(r,X(r,x^*))|.$ Then by Lemma \ref{rev2}, $\Phi_{g} \in {\bf D}_p^1.$ Let $\Phi_{g,i}(r,x), i = 1, \ldots, N$ be defined by the formula of Lemma \ref{rev}.
Then we have
$$I_{1,j}( r,F) = r^{-1/2} E^{\mu} [ \Phi_{g,j}(r,x) \varphi_m \circ ( P_{T-r}F)(X(r,x^*))  V_j (P_{T-r}F) (X(r,x^*))],$$
and 
$$
\sup_{t \in[0,T], x\in {\bf R}^N} E^{\mu}[|\Phi_{g,i}(t,x)|^p] < \infty.
$$ 
Then there exists a constant $C>0$ such that
$$|I_{1,j}( r,F)| \leqq C r^{-1/2} \|V_j (P_{T-r}F) \|_{\infty}.$$
Also we have
$$|I_{2,j}( r,F)| \leqq  CE^{\mu}[|g(r, X(r,x^*))|]\|V_j^2 (P_{T-r}F)\|_{\infty},$$
for any $ F \in C_b^{\infty}({\bf R}^{N})$ and any $0< r< T$.\\ 
\\Let vector field $V_j$ be represented by   
$ V_j = \sum_{i=1}^N v_j^i(x) \frac{\partial}{\partial x_i}.$
Then we have $$V_j (P_{T-r} F)(x)= \sum_{i=1}^N\sum_{k=1}^N v_j^i(x) (T_{\Phi_{k,i}}(T-r)\frac{\partial F}{\partial x_i})(x),$$
where $\Phi_{k,i}(t,x)=\frac{\partial X^{k}(t,x)}{\partial x_i}$ and 
$$(T_{\Phi_{k,i}}(t)F)(x) = E^{\mu}[\Phi_{k,i}(t,x) F(X(t,x))]. $$
Moreover, we have  
$$V_j^2 (P_{T-r} F)(x)= \sum_{i=1}^N\sum_{k=1}^N ( V_j v_j^i(x)  (T_{\Phi_{k,i}}(T-r)\frac{\partial F}{\partial x_i})(x)
+ v_j^i(x)  (V_j T_{\Phi_{k,i}}(T-r)\frac{\partial F}{\partial x_i})(x)).$$
Then by Corollary 9 of \cite{K}, since $\Phi_{k,i} \in {\cal K}_0$ and there is a constant $C>0$ such that
$$\|V_j (P_{T-r} F) \|_{\infty} \leqq C\|\nabla F\|_{\infty},$$
and
$$\|V_j^2 (P_{T-r} F)  \|_{\infty} \leqq C(T-r)^{-1/2}\|\nabla F\|_{\infty},$$
for any $ F \in C_b^{\infty}({\bf R}^{N}), j = 1,\ldots,d,$ and any $0< r< T$.\\ 
\\
So we have
$$E[|g(t,X(t,x^*))\bar{\varphi}_m((P_{T-t}F)(X(t,x)) ) - g(s,X(s,x^*))\bar{\varphi}_m((P_{T-s}F)(X(s,x)))|] $$
$$\leqq C \|\nabla F\|_{\infty} \int_s^t (r^{-1/2} + (T-r)^{-1/2})dr.$$
Letting $m \to \infty,$ we have our assertion.
 \qed
 
\begin{cor}
Let $T>0$. There exists a constant $C >0$ such that
$$E[|g(t,X(t,x^*))(P_{T-t}F)(X(t,x^*)) \vee 0 - g(s,X(s,x^*))(P_{T-s}F)(X(s,x^*)) \vee 0|]$$
$$\leqq C \|F\|_{{Lip}} \int_{s}^t(r^{-1/2}+(T-r)^{-1/2})dr ,$$
for any $ F \in Lip({\bf R}^{N})$ and any $0< s < t< T$.
\end{cor}
{\it Proof.}
For $F \in Lip({\bf R}^N)$, there exists $F_m \in C_b^{\infty}({\bf R}^N), m=1, 2, \cdots, $ such that
$\| \nabla F_m\| \leqq \|F\|_{Lip}$ and $F_m(x) \to F(x)$, for any $x \in {\bf R}^N$.
So we obtain the result from Lemma \ref{discLem2}.
\qed
\\
\begin{lem}\label{discLem1}
Let $m=1,\ldots,M,$ and $T>0.$ There exists a constant $C >0$ such that
\begin{align}
&E^{\mu}[ |g(t,X(t,x^*))(P^{(m)}_{T-t}h)(\tilde{X}^{(m)}(t,\tilde{x}_m^*))-h(\tilde{X}^{(m)}(t,\tilde{x}_m^*))|]\nonumber \\ 
\leqq &C(\|\nabla h\|_{\infty}+ \|\nabla^2 h\|_{\infty} )(T-t),\label{disc1}\\
\text{and} \nonumber\\
&E^{\mu}[ |g(t,X(t,x))(P^{(m)}_{T-t}(h\vee 0))(\tilde{X}^{(m)}(t,\tilde{x}_m^*))-(h\vee 0)(\tilde{X}^{(m)}(t,\tilde{x}_m^*))|] \nonumber \\
\leqq &C(\|\nabla h\|_{\infty}+ \|\nabla^2 h\|_{\infty} )(T-t) \label{disc2}.
\end{align}
for  any $h \in C_b^{\infty}({\bf R}^{\tilde{N}_m}), t \in [0, T).$
\end{lem}
{\it Proof.}
(\ref{disc1}) follows from It\^o's formula. So we show (\ref{disc2}).\\
Let $\bar{\varphi}_k, k=1,\ldots,$ are as defined in (\ref{defPhi}). Let
$$\tilde{L}_m=\tilde{V}^{(m)}_{0}+\frac{1}{2}\sum_{i=1}^{d}(\tilde{V}^{(m)}_{i})^2.$$
By It\^o's formula
$$\bar{\varphi}_k(h(\tilde{X}^{(m)}(T,\tilde{x}_m^*))-\bar{\varphi}_k(h(\tilde{X}^{(m)}(t,\tilde{x}_m^*))$$
$$=\int_t^T \varphi_k(h(\tilde{X}^{(m)}(s,\tilde{x}_m^*)) (\tilde{V}^{(m)}_{i}h)(\tilde{X}^{(m)}(s,\tilde{x}_m^*)) dB^{m,i}(s)$$
$$+\int_t^T \varphi_k(h(\tilde{X}^{(m)}(s,\tilde{x}_m^*)) (\tilde{L}_m h)(\tilde{X}^{(m)}(s,\tilde{x}_m^*)) ds$$
$$+\frac{1}{2} \int_t^T \varphi_k' (h(\tilde{X}^{(m)}(s,\tilde{x}_m^*)) \sum_{i=1}^{\tilde{d}_m} ((\tilde{V}^{(m)}_{i} h)(\tilde{X}^{(m)}(s,\tilde{x}_m^*)))^2 ds$$

So we have
$$E^{\mu}[\bar{\varphi}_k(h(\tilde{X}^{(m)}(T,\tilde{x}_m^*))|\tilde{\mathcal{F}}^{(m)}_t]-\bar{\varphi}_k(h(\tilde{X}^{(m)}(t,\tilde{x}_m^*))$$
$$=\int_t^T E^{\mu}[\varphi_k(h(\tilde{X}^{(m)}(s,\tilde{x}_m^*)) (\tilde{L}_m h)(\tilde{X}^{(m)}(s,\tilde{x}_m^*))|\tilde{\mathcal{F}}^{(m)}_t] ds$$
$$+\frac{1}{2} \sum_{i=1}^{\tilde{d}_m}  \int_t^T E[\varphi_k' (h(\tilde{X}^{(m)}(s,\tilde{x}_m^*))) ((\tilde{V}^{(m)}_{i} h)(\tilde{X}^{(m)}(s,\tilde{x}_m^*)))^2 |\tilde{\mathcal{F}}^{(m)}_t] ds$$
Notice that $\varphi_k' \geqq 0, $ then we have
$$E[|g(t,X(t,x^*))E^{\mu}[\bar{\varphi}_k(h(\tilde{X}^{(m)}(T,\tilde{x}_m^*))|\tilde{\mathcal{F}}^{(m)}_t]-\bar{\varphi}_k(h(\tilde{X}^{(m)}(t,\tilde{x}_m^*))|]$$
$$\leqq E[|g(t,X(t,x^*))|]\|\tilde{L}_m h\|_{\infty}(T-t)$$
$$+\frac{1}{2} \sum_{i=1}^{\tilde{d}_m} \int_t^T E[|g(t,X(t,x^*))| \varphi_k' (h(\tilde{X}^{(m)}(s,\tilde{x}_m^*)) (\tilde{V}^{(m)}_{i} h(\tilde{X}^{(m)}(s,\tilde{x}_m^*)))^2] ds.$$
On the other hand, we have
$$ \varphi_k' (h(\tilde{x}_m^*)) (\tilde{V}^{(m)}_{i} h(\tilde{x}_m^*))^2$$
$$= \tilde{V}^{(m)}_{i} (\varphi_k \circ h)  (\tilde{x}_m^*) \left(\tilde{V}^{(m)}_{i} h\right)(\tilde{x}_m^*) $$
$$= \tilde{V}^{(m)}_{i} \left( (\varphi_k \circ h) (\tilde{x}_m^*)  \tilde{V}^{(m)}_{i} h(\tilde{x}_m^*) \right) 
-(\varphi_k \circ h) (\tilde{x}_m^*)  (\tilde{V}^{(m)}_{i})^2 h (\tilde{x}_m^*) . $$ 
\\
Let $\Phi_g(t,x) = |g(t,X(t,x^*))|$ and $\Phi_{g,i}(t,x), i = 1,\ldots,N$ be defined by the formula of Lemma \ref{rev}.
Then it follows that
$$ \|E^{\mu} [\Phi_g(t,x) \tilde{V}^{(m)}_{i} \left( (\varphi_k \circ h(\tilde{X}^{(m)}(s,\tilde{x}_m^*)))) \tilde{V}^{(m)}_{i} h(\tilde{X}^{(m)}(s,\tilde{x}_m^*))) \right) ]\|_{\infty} $$
$$\leqq Cs^{-1/2} E^{\mu}[|\Phi_{g,1}(t,x)|] \|(\varphi_k \circ h) \tilde{V}^{(m)}_{i} h \|_{\infty} \leqq C s^{-1/2} \| \tilde{V}^{(m)}_{i} h \|_{\infty},$$
and
$$ \| E^{\mu}[ \Phi_g(t,x) (\varphi_k \circ h) (\tilde{X}^{(m)}(s,\tilde{x}_m^*))  (\tilde{V}^{(m)}_{i})^2 h (\tilde{X}^{(m)}(s,\tilde{x}_m^*)) ] \|_{\infty} 
\leqq C\| (\tilde{V}^{(m)}_{i})^2 h \|_{\infty}.$$
So we have
$$\frac{1}{2} \sum_{i=1}^{\tilde{d}_m} \int_t^T E[g(t,X(t,x))\varphi_k' (h(\tilde{X}^{(m)}(s,\tilde{x}_m^*)) (\tilde{V}^{(m)}_{i} h(\tilde{X}^{(m)}(s,\tilde{x}_m^*)))^2] ds$$
$$\leqq  \int_t^T C' (\|\nabla h\|_{\infty}+ \|\nabla^2 h\|_{\infty} )(1+s^{-1/2}) ds $$
$$\leqq C' (\|\nabla h\|_{\infty}+ \|\nabla^2 h\|_{\infty} )(T-t) (1+(T^{1/2}+t^{1/2})^{-1}).$$
Letting $k \to \infty,$ we have the assertion. \qed

\begin{cor}
Let $m=1,\ldots,M,T>0$ and $F \in \mathcal{M}({\bf R}^{\tilde{N}_m})$. There exists a constant $C>0$ such that
\begin{align}
&E^{\mu}[ |g(t,X(t,x^*))(P^{(m)}_{T-t}F)(\pi_m{X}(t,x^*))-F(\pi_m{X}(t,x^*))|] 
\leqq C(T-t),
\end{align}  
for any $ t \in [0,T).$
\end{cor}
{\it Proof.}
Notice that $\pi_m{X}(t,x)=\tilde{X}^{(m)}(t,\tilde{x}_m^*)$ and Lemma \ref{discLem1} is valid for $h \in \hat{\mathcal{D}}({\bf R}^{\tilde{N}_m})$.
On the other hand, for $F\in \mathcal{M} ({\bf R}^{\tilde{N}_m}) $, 
we have the expression that 
$$F=\sum_{k=1}^{K_F} a_k (f_k\vee g_k) = \sum_{k=1}^{K_F} a_k\left( (f_k-g_k)\vee 0 +g_k\right),$$
$a_k \in {\bf R}, f_k, g_k \in  \hat{\mathcal{D}}({\bf R}^{\tilde{N}_m}), k=1,\ldots, K_F.$
So our assertion follows from Lemma \ref{discLem1}.
 \qed

\section{Discretization}
Let $c_{\Delta}, \Delta \in \Pi,$ be given by
$$c_{\Delta}=\sum_{i=0}^{n-1}(t_{i+1}-t_i)
E^{\mu}[g(t_i,X(t_i,x^*))\{\sum_{m=1}^M\sum_{k; T_k \geqq t_{i+1}}^K(P_{T_k-t_i}^{(m)}F_{m,k})(\pi_mX(t_i,x^*))) \}\vee 0].
$$
Let $i_{(k)}, k=1,\ldots, K,$ be such that $T_k=t_{i_{(k)}}.$
Then we have $c_{\Delta}$ is as follows.
$$
c_{\Delta}=\sum_{k=1}^K\sum_{i=i_{(k-1)}}^{{i_{(k)}}-1}(t_{i+1}-t_i) E^{\mu}[ g(t_i,X(t_i,x^*))\{
\sum_{m=1}^M \sum_{k'=k}^K(P_{T_k-t_i}^{(m)}F_{m, k'})(\pi_{m} X(t_i,x^*))\}\vee0 ]
 $$

Let ${\cal F}_{t}^{(\infty)},$ $t\geqq 0,$  be 
sub $\sigma$-algebra of ${\cal F}$ given by
$$
{\cal F}_{t}^{(\infty )} 
= \sigma \{ X_{\ell}(s); \; s\in [0,t], \ \ell =1,2,\ldots \} .
$$
\\
\begin{prop}\label{time disc}
There exists a constant $C >0$ such that
$$
|c_0-c_{\Delta}| \leqq  C |\Delta|,\quad \Delta \in \Pi$$.
\end{prop}
{\it Proof.}
Let $$\tilde{F}_{k}(x)=\sum_{m=1}^M \sum_{k'=k}^K (P_{T_{k'}-T_{k}}^{(m)}F_{m,k'})(\pi_m x), k=1,\ldots,K.$$ 
Then by Lemma \ref{discLem2}, there is a constant $C>0$ such that\\
\begin{align*}
&|c_0- c_{\Delta}| \\
 \leqq   &\sum_{k=1}^K \sum_{i=i_{(k-1)}}^{i_{(k)}-1} \int_{t_i}^{t_{i+1}} |E^{\mu}[ g(t,X(t,x^*))\left( (P_{T_{k}-t} \tilde{F}_{k})(X(t,x^*)) \vee 0\right)] \\
-& E^{\mu}[ g(t_i,X(t_i,x^*))(P_{T_{k}-t_i} \tilde{F}_{k})(X(t_i,x^*)) \vee 0]| dt\\
\leqq   &C\sum_{k=1}^K \sum_{i=i_{(k-1)}}^{i_{(k)}-1} \int_{t_i}^{t_{i+1}} dt \int_{t_i}^t (r^{-1/2}+(T_k-r)^{-1/2})dr\\
\leqq &C|\Delta|  \sum_{k=1}^K \sum_{i=i_{(k-1)}}^{i_{(k)}-1}\int_{t_i}^{t_{i+1}} (r^{-1/2}+(T_k-r)^{-1/2})dr. 
\end{align*}
So we have
$$
|c_0- c_{\Delta}| 
\leqq 
C|\Delta| \sum_{k=1}^K \int_{T_{k-1}}^{T_k} (r^{-1/2}+(T_k-r)^{-1/2})dr .
$$
So the assertion follows.
\qed

\section{Property of Stochastic Mesh Operator}

To estimate the stochastic mesh operator, we use the following estimation of 
transition kernel $p^{(m)}(t,\tilde{x}_m, x)$ obtained by Proposition 8 of  \cite{KM}.
\begin{proposition}\label{heatKernel}
Let $\delta_0^{(m)}$ be given by 
$$\delta_0^{(m)}=(3\tilde{N}_m(\sup_{x\in {\bf R}^{\tilde{N}_m}} \sum_{i=1}^{d}|\tilde{V}^{(m)}_{i}(x)|^2))^{-1},$$
then for any $T>0,$ and $m=1,\ldots,M,$ there is a $C>0$ such that
$$
p^{(m)}(t,x,y)
\leqq C t^{- (\tilde{N}_m+1) \ell_0 /2} h_m(x)^{-2(\tilde{N}_m+1)\ell_0}\exp (- \frac{2\delta_0^{(m)}}{t}|y-x|^2 ), 
\qquad  t \in (0,T], \ x,y \in E_m,
$$
and
$$
p^{(m)}(t,x,y)
\leqq C t^{- (\tilde{N}_m+1) \ell_0 /2} h_m(y)^{-2(\tilde{N}_m+1)\ell_0}\exp (- \frac{2\delta_0^{(m)}}{t}|y-x|^2 ), 
\qquad  t \in (0,T], \ x,y \in E_m.
$$
In particular, for any $T>0, m=1,\ldots,M,$
 and $q\geqq 1,$ 
there is a $C>0$ such that
$$
p^{(m)}(t,x,y)
\leqq C t^{- (\tilde{N}_m+1) \ell_0/2} h_m(x)^{-2(\tilde{N}_m+1)\ell_0}(1+|x|^2)^q
(1+|y|^2)^{-q}, 
\qquad  t \in (0,T], \ x,y \in E_m.
$$
\end{proposition}

Let $\nu^{(m)}_t(dx)=p^{(m)}(t, \tilde{x}^*_m, x)dx.$
From Proposition 13, 21 and Proposition 15 (1) of \cite{KM}, we have the followings.
\begin{proposition}\label{basic}
Let $t > 0, f\in L^2(E_m;d\nu_t^{(m)})$ and $t>s\geqq 0.$ Then we have
$$
E^P[(Q_{s,t}^{(m,L,\omega)}f)(x)| {\cal F}_{s}^{(\infty )}] 
= (P^{(m)}_{s,t}f)(x),
\qquad \nu_s^{(m)}-a.e.x \in E_m.
$$
and 
$$
E^P[ |(Q_{s,t}^{(m,L,\omega)}f)(x) - (P^{(m)}_{s,t}f)(x)|^2 | {\cal F}_s^{(\infty )}]
\leqq 
\frac{1}{L}
\int_{E_m} \frac{p^{(m)}(t-s,x,y)^2 |f(y)|^2}{q_{s,t}^{(m,L,\omega)}(y)}\; dy .
$$
\end{proposition}

\begin{proposition}\label{key}
Let $\delta \in (0,1)$ then there exists a $C > 0$ such that
\begin{align*}
&\left( \frac{1}{L}\sum_{\ell=1}^L E^P[ |(Q_{t,T_k,\varepsilon}^{(m)}f)(\pi_m (X_{\ell}(t))) - (P^{(m)}_{t,T_k}f)(\pi_m (X_{\ell}(t))|^2 ]\right)^{1/2} \\
&\leqq 
C  L^{-(1-\delta)/2}(T_k-t)^{-(1+\delta)(\tilde{N}+1)\ell_0/4} (\int_{E_m} f(y)^2 (1+|y|^2)^{-\tilde{N}_m}dy)^{1/2}.
\end{align*}
for any $\varepsilon >0$ any $m=1,\ldots,M,$  and any $f \in Lip({\bf R}^{\tilde{N}_m}).$
\end{proposition}

\begin{prop}\label{exSet}
Let
\begin{align*}
Z_L^{(m,k)}(t;\delta)=\sup_{y\in {\bf R}^{\tilde{N}_m}} \frac{|q_{t,T_k}^{(m,L,\omega)}(y)- p^{(m)}(T_k, \tilde{x}^*, y)|}
{(L^{-1/(1-\delta)}+p^{(m)}(T_k, \tilde{x}^*,y))^{(1-\delta)/2}},
\end{align*}
\begin{align*}
\tilde{Z}_L^{(m,k)}(t; \delta)=\sup_{s\in [0,t]} Z_L^{(m,k)}(s;\delta).
\end{align*}
Then we have the followings.\\
$(1)$ For any $\delta \in (0,1)$, and $p>1$ , there is a $C_{p,\delta} > 0$ such that
$$
E^P[(L^{(1-\delta^2)/2}\tilde{Z}_L^{(m,k)}(T_k-\varepsilon; \delta))^p]^{1/p}\leqq C_{p,\delta}\varepsilon^{-5\ell_0}L^{-p\delta^2/2+1/p}$$
for any $ \varepsilon \in  (0,T_k], k=1,\ldots,K,$  and  $L \geqq 1.$\\
$(2)$ Let $\delta \in (0,1), t \in (0, T_k)$ and $\varepsilon \in (0,T).$  If $L^{(1-\delta^2)/2}\tilde{Z}_L^{(m,k)}(t; \delta) \leqq 1/4$, and  $p^{(m)}(T_k, x_0, y) \geqq L^{-(1-\delta)}$, then
$$
\frac{1}{2} \leqq \frac{q_{t,T_k}^{(m,L,\omega)}(y)}{p^{(m)}(T_k, \tilde{x}^*, y)} \leqq 2, $$
for any $t \in (0, T_k-\varepsilon], \quad k=1,\ldots,K.$ and $L \geqq 1.$
\end{prop}
Now we introduce the following sets and functions.
Let $B^{(m,k)}(t, \delta, L) \in \mathcal{F},$ $\varphi_{m,k,L}, m=1,\ldots,M, k=1, \ldots,K,$ be given by
$$
B^{(m,k)}(t, \delta,L) = \{\omega \in \Omega; L^{(1-\delta^2)/2}\tilde{Z}_L^{(m,k)}(t; \delta)
\leqq 1/4 \},
$$
and
$$
\varphi_{m,k,L}(y) = 1_{\{ y \in E_m; p^{(m)}(T_k, x_0, y) > L^{-(1-\delta)} \}}.
$$
Let $d_{i,\varepsilon,L}^{(m,k)}: [0,T] \times E \times \Omega \to [0, \infty), i=1,2,3,$ be the measurable functions given by
\begin{align*}
 d_{1,\varepsilon,L}^{(m,k)}(t,x)& = |(Q_{t,T_k, \varepsilon}^{(m,L,\omega)}(1-\varphi_{m,k,L})F_{m,k})(\pi_m(x))
-(P_{T_k-t}^{(m)} (1-\varphi_{m,k,L})F_{m,k})(\pi_m(x)) |1_{[ 0, T_k- \varepsilon ) }(t),\\
 d_{2,\varepsilon,L}^{(m,k)}(t,x) &= 1_{B^{(m,k)}(T_k-\varepsilon, \delta,L)}|(Q_{t,T_k, \varepsilon}^{(m,L,\omega)}\varphi_kF_{m,k})(\pi_m(x))
-(P_{T_k-t}^{(m)} \varphi_{m,k,L}F_{m,k})(\pi_m(x)) |1_{[ 0, T_k- \varepsilon ) }(t),\\
 d_{3,\varepsilon,L}^{(m,k)}(t,x) &= |(Q_{t,T_k, \varepsilon}^{(m,L,\omega)}F_{m,k})(\pi_m(x))
-(P_{T_k-t}^{(m)} F_{m,k})(\pi_m(x)) |1_{[T_k - \varepsilon , T_k) }(t)\\
&= |F_{m,k}(\pi_m(X(T_k,x^*)))
-(P_{T_k-t}^{(m)} F_{m,k})(\pi_m(x)) |1_{[T_k - \varepsilon , T_k) }(t), \quad k = 1, \cdots, K.
\end{align*}
 Let $p(t,x,dy)$ be the transition kernel of $X(t,x).$
\begin{prop} \label{d12}
Let $\delta \in (0,1)$. Then there exists a constant $C > 0$ such that 
\begin{align}
&\int_{E} E^P[ d_{1, \varepsilon,L}^{(m,k)}(t,x) ]|g(t,x)| p(t, x^*, dx)
\leqq 
 C L^{-(1-\delta)^3} 1_{[ 0, T_k- \varepsilon ) }(t), \label{d1} \\
& (\int_{E_m} E^P[\ d_{2,\varepsilon,L}^{(m,k)}(t,x)^2] p(t, x^*, dx))^{1/2}\nonumber \\
\leqq  
&C L^{-(1-\delta)/2} (T_k-t)^{-(1+\delta)(\tilde{N}+1)\ell_0/4} 1_{[ 0, T_k- \varepsilon ) }(t),\label{d2}\\
\text{and} \nonumber\\
&  \int_{E_m}  d_{3,\varepsilon,L}^{(m,k)}(t,x) |g(t,x)| p(t, x^*, dx) 
\leqq C (T_k-t)1_{[T_k - \varepsilon , T_k) }(t). \label{d3}
\end{align}
for any $\varepsilon \in (0,\varepsilon_0), t \in (0, T_k], L \geqq 1,$ $m=1,\ldots,M,$ and $k=1,\ldots,K.$
\end{prop}
{\it Proof.}
Equation (\ref{d3}) follows from Lemma \ref{discLem1}. 
So we will show  (\ref{d1}) and (\ref{d2}). Note that if $t \geqq T_k - \varepsilon$, 
both sides are $0$ in (\ref{d1}) and (\ref{d2}) . So we will consider the case $t < T_k - \varepsilon$.
By Proposition \ref{basic}, we have
\begin{align*}
&\int_{E} E^P[d_{1,\varepsilon,L}^{(m,k)}(t,\pi_m(x))] |g(t,x)| p(t, x^*, dx)\\
= &\int_{E} E^P[|(Q_{t,T_k, \varepsilon}^{(m)}(1-\varphi_{m,k,L})F_{m,k})(x)
-(P_{T_k-t}^{(m)} (1-\varphi_{m,k,L})F_{m,k})(\pi_m(x)) |] |g(t,x)| p(t, x^*, dx) \\
\leqq &\int_{E} (E^P[(Q_{t,T_k, \varepsilon}^{(m)} |(1-\varphi_{m,k,L})F_{m,k}| )(\pi_m(x))]\\
& \qquad \qquad \qquad \qquad \qquad \qquad+ (P_{T_k-t}^{(m)} |(1-\varphi_{m,k,L})F_{m,k}| )(\pi_m(x))) |g(t,x)| p(t, x^*, dx)\\
\leqq &2 \int_{E} (P_{T_k-t}^{(m)} (1-\varphi_{m,k,L})|F_{m,k}|)(\pi_m(x)) g(t,x) p(t, x^*, dx).
\end{align*}
Using H\"{o}lder's inequality for $p=\frac{1}{\delta}, q = \frac{1}{1-\delta},$
$$\int_{E} (P_{T_k-t}^{(m)} (1-\varphi_{m,k,L})|F_{m,k}|)(\pi_m(x)) |g(t,x)| p(t, x^*, dx)$$
$$\leqq  \{\int_{E} (P_{T_k-t}^{(m)} (1-\varphi_{m,k,L})|F_{m,k}|)(\pi_m(x))^{1/(1-\delta)} p(t, x^*, dx)\}^{1-\delta}$$
$$\times \{\int_E |g(t,x)|^{1/\delta}p(t, x^*, dx)\}^{\delta}$$
$$\leqq  \{\int_{E_m}  (1-\varphi_{m,k,L}(\tilde{y}_m))^{1/(1-\delta)}|F_{m,k}(\tilde{y}_m)|^{1/(1-\delta)} p^{(m)}(T_k, \pi_m(x^*), \tilde{y}_m)d\tilde{y}_m\}^{1-\delta}$$
$$\times \{\int_E |g(t,x)|^{1/\delta}p(t, x^*, dx)\}^{\delta}$$
$$\leqq  L^{-(1-\delta)^3}\{  \int_{E_m} |F_{m,k}(\pi_m(y))|^{1/(1-\delta)} p^{(m)}(T_k, \pi_m(x^*),\tilde{y}_m)^{\delta}d\tilde{y}_m \}^{1-\delta}$$
$$\times \{\int_E |g(t,x)|^{1/\delta}p(t, x^*, dx)\}^{\delta}.$$
We used $(1-\varphi_{m,k,L}(\tilde{y}_m))^{1/(1-\delta)} p^{(m)}(T_k, \pi_m(x^*), \tilde{y}_m)^{(1-\delta)} \leqq L^{-(1-\delta)^2}$ in the last inequality.
So we have Equation (\ref{d1}).\\
Next we will show Equation (\ref{d2}) . 
Noting that from
$B^{(m,k)}(T_k-\varepsilon, \delta,L) \subset B^{(m,k)}(t, \delta,L), t \in [0,T_k-\varepsilon), k=1,\ldots,K,$ and $L \geqq 1$,
$$d_{2,\varepsilon,L}^{(m,k)}(t,x) 
\leqq 1_{B^{(m,k)}(t, \delta,L)}|(Q_{t,T_k, \varepsilon}^{(m)}\varphi_{m,k,L}F_{m,k})(\pi_m(x))
-(P_{T_k-t}^{(m)} \varphi_{m,k,L}F_{m,k})(\pi_m(x))|.$$ 
Since Propsition \ref{exSet},  
$ 1_{B^{(m,k)}(t, \delta,L)} q_{t,T_k}^{(m,L,\omega)}(\tilde{y}_m)^{-1} \leqq 2 p^{(m)}(T_k, \pi_m(x^*),\tilde{y}_m)^{-1}$. And
by Proposition \ref{basic}, 
we have
\begin{align*}
&1_{B^{(m,k)}(t, \delta,L)}E^P[ |(Q_{t,T_k, \varepsilon}^{(m)}\varphi_{m,k,L}F_{m,k})(\tilde{x}_m)
-(P_{T_k-t}^{(m)} \varphi_{m,k,L}F_{m,k})(\tilde{x}_m) |^2| \mathcal{F}_t]\\
\leqq &1_{B^{(m,k)}(t, \delta,L)} \frac{1}{L} \int_{E_m} \frac{|\varphi_{m,k,L}F_{m,k}(\tilde{y}_m)|^2p^{(m)}(T_k-t, \tilde{x}_m,\tilde{y}_m)^2}{q_{t,T_k}^{(m,L,\omega)}(\tilde{y}_m)}d\tilde{y}_m\\
\leqq &\frac{2}{L} \int_{E_m} \frac{|\varphi_{m,k,L}F_{m,k}(\tilde{y}_m)|^2}{p^{(m)}(T_k, \pi_m(x^*),\tilde{y}_m)} p^{(m)}(T_k-t, \tilde{x}_m,\tilde{y}_m)^2 dy.
\end{align*}
Then we have
\begin{align*}
&(\int_{E} E^P[\ d_{2,\varepsilon,L}^{(m,k)}(t,x)^2] p(t, x^*, dx))^{1/2} \\
\leqq &(\int_{E}  E^P[1_{B^{(m,k)}(t, \delta,L)}E^P[ |(Q_{t,T_k, \varepsilon}^{(m)}\varphi_{m,k,L}F_{m,k})(x)
-(P_{T_k-t}^{(m)} \varphi_{m,k,L}F_{m,k})(x) |^2| \mathcal{F}_t] ]p(t, x^*, dx))^{1/2}\\
\leqq &(\frac{2}{L} \int_{E_m} \frac{|\varphi_{m,k,L}F_{m,k}(y)|^2}{p^{(m)}(T_k, \pi_m(x^*),\tilde{y}_m)} 
(\int_{E} p^{(m)}(T_k-t, \tilde{x}_m,\tilde{y}_m)^{(1-\delta)+(1+\delta)} p(t, x^*, dx))d\tilde{y}_m)^{1/2}\\
\leqq &(\frac{2}{L} \int_{E_m} \frac{|\varphi_{m,k,L}F_{m,k}(y)|^2}{p^{(m)}(T_k, \pi_m(x^*),\tilde{y}_m)^{\delta}} 
(\int_{E_m} p^{(m)}(t, \pi_m(x^*), \tilde{x}_m) p^{(m)}(T_k-t, \tilde{x}_m,\tilde{y}_m)^{(1+\delta)/\delta}d\tilde{x}_m)^{\delta}d\tilde{y}_m)^{1/2}.
\end{align*}
Let $q \geqq \tilde{N}.$ From Lemma \ref{heatKernel}, there exists a constant $C>0$ such that 
$$p^{(m)}(T_k-t, \tilde{x}_m,\tilde{y}_m) \leqq C(T_k-t)^{-(\tilde{N}_m+1)\ell_0/2}h_m(\tilde{x}_m)^{-(\tilde{N}_m+1)\ell_0}(1+|\tilde{x}_m|^2)^q(1+|\tilde{y}_m|^2)^{-q}.$$
We set $C_1$ as
\begin{align*}
C_1=&\sup_{t\in[0,T]}\max_{\substack{m=1,\ldots,M,\\ k=1,\ldots,K}}(\int_{E_m} h_m(x)^{-(\tilde{N}_m+1)\ell_0(1+\delta)/\delta}(1+|x|^2)^{q(1+\delta)/\delta}p^{(m)}(t, \tilde{x}_m^*,x)dx)^{\delta/2}\\
&\quad \quad \times (\int_E|g(t,x)|dx)^{1/2}.
\end{align*}
$C_1$ is bounded by Proposition 3 of \cite{KM}.
Then since $\varphi_{m,k,L}(y)p^{(m)}(T_k, \tilde{x}_m^*,y)^{-\delta} \leqq L^{\delta}$,  we have
\begin{align*}
&\int_{E} E^P[\ d_{2,\varepsilon,L}^{(m,k)}(t,x)^2]^{1/2}|g(t,x)| p(t, x^*, dx)\\
\leqq &\frac{C_1}{L} \int_{E_m} p^{(m)}(T_k, \tilde{x}_m^*,\tilde{y}_m)^{-\delta} |\varphi_{m,k,L}F_{m,k}(\tilde{y}_m)|^2(1+|\tilde{y}_m|^2)^{-q(1+\delta)}d\tilde{y}_m 
(T_k-t)^{-(1+\delta)(N+1)\ell_0/2} ,\\
\leqq &C_1L^{-(1-\delta)}(T_k-t)^{-(1+\delta)(\tilde{N}_m+1)l_0/2} 
\int_{E_m} |F_{m,k}(\tilde{y}_m)|^2(1+|\tilde{y}_m|^2)^{-q(1+\delta)}d\tilde{y}_m .
\end{align*}
Since $q \geqq \tilde{N},$ and $F_{m,k}$ is Lipschitz continuous,
$$\int_{{\bf R}^{\tilde{N}_m}} |F_{m,k}(\tilde{y}_m)|^2(1+|\tilde{y}_m|^2)^{-q(1+\delta)}d\tilde{y}_m < \infty.$$
So we have the assertion.
\qed
\\
Let $a,b,\alpha,\beta \geqq 0,$ and $a_i,b_i,\alpha_i,\beta_i \geqq 0, i=1,2.$
Let $\phi^{(k)}(t,\varepsilon; a, \alpha, b, \beta)$ and $\hat{e}(\varepsilon, \gamma), t \in [0, T_k)$ be
\begin{align*}
&\phi^{(k)}(t,\varepsilon; a, \alpha, b, \beta) = a(T_k-t)^{-\alpha} 1_{[0, T_k - \varepsilon ) }(t)
+b(T_k-t)^{\beta}1_{[T_k - \varepsilon , T_k) }(t),\\
&\hat{e}(\varepsilon, \gamma)
=
\begin{cases}
 \varepsilon^{-(\gamma-1)} )&, \gamma >1, \\
 \log \varepsilon &, \gamma = 1, \\
 1&, 0\leqq \gamma < 1.\\
\end{cases}
\end{align*}

\begin{prop}\label{dd} 
There exists a constant $C > 0$ such that
\begin{align} \label{sum1}
 \sum_{i=0}^{n-1} (t_{i+1}-t_i) \sum_{k; T_k \geqq t_{i+1}} \phi^{(k)}(t_i, \varepsilon; a, \alpha, b, \beta)
\leqq C (a\hat{e}(\varepsilon,\alpha ) + b \varepsilon^{\beta+1}),
\end{align}
and,
\begin{align} \label{sum2}
 &\sum_{i=0}^{n-1} (t_{i+1}-t_i) (\sum_{k_1; T_{k_1} \geqq t_{i+1}}  \phi^{(k)}(t_i, \varepsilon; a_1, \alpha_1, b_1, \beta_1)) 
 (\sum_{k_1; T_{k_1} \geqq t_{i+1}} \phi^{(k)}(t_i, \varepsilon; a_2, \alpha_2, b_2, \beta_2)) \nonumber \\
\leqq &C \left( a_1 a_2 \hat{e}(\varepsilon,\alpha_1+\alpha_2 ) +
a_1 b_2 \varepsilon^{\beta2+1} 
+a_2 b_1 \varepsilon^{\beta_1+1} 
+ b_1 b_2 \varepsilon^{\beta_1+ \beta2+1} \right)
\end{align}
for any $\varepsilon >0$.
\end{prop}

\noindent {\it Proof.}
Let us take $i_{(k)}$ as $t_{i_{(k)}}=T_k, k=1,\ldots,K.$
If $t_i \in [T_{k-1}, T_k]$ and $k'>k$ then $ T_{k'}- t_i > \varepsilon .$
So notice that
\begin{equation}\label{discEps}
1_{[T_{k'}-\varepsilon , T_{k'})}(t_i) =0,\quad  \text{for} \ i_{(k-1)} \leqq i \leqq i_{(k)}-1, k' >k.
\end{equation}
So we have
\begin{align*}
&\sum_{i=0}^{n-1} (t_{i+1}-t_i) \sum_{k; T_k \geqq t_{i+1}} \phi^{(k)}(t_i, \varepsilon; a, \alpha, b, \beta)
= \sum_{k=1}^K \sum_{i=i_{(k-1)}}^{i_{(k)}-1} (t_{i+1}-t_i) \sum_{k'=k}^K \phi^{(k')}(t_i, \varepsilon; a, \alpha, b, \beta ) \\
=&\sum_{k=1}^K \sum_{i=i_{(k-1)}}^{i_{(k)}-1} (t_{i+1}-t_i) 
\sum_{k'=k}^K  \left( a(T_{k'}-t_i)^{-\alpha} 1_{[0, T_k - \varepsilon ) }(t_i)
+b (T_{k'}-t_i)^{\beta}1_{[T_{k'} - \varepsilon , T_{k'}) }(t_i) \right)\\
\leqq &\sum_{k=1}^K  \left( \sum_{k'=k}^K \int_{T_{k'-1}}^{T_{k'}-\varepsilon} a(T_{k'}-t)^{-\alpha} dt
+ \sum_{i=i_{(k-1)}}^{i_{(k)}-1} b1_{ [T_k- \varepsilon ,  T_k)}(t_i) (t_{i+1}-t_i)(T_k-t_i)^{\beta} \right), 
\end{align*}
because $(T_{k'}-t_i)^{-\alpha} \leqq (T_{k'}-t)^{-\alpha}$ for  $t_i \leqq t \leqq t_{i+1} $.
\\
On the other hand,
$$ 
\sum_{k=1}^K   \sum_{k'=k}^K \int_{T_{k'-1}}^{T_{k'}-\varepsilon} (T_{k'}-t)^{-\alpha} dt 
\leqq K^2  \hat{e}(\varepsilon,\alpha),$$
and
$$
\sum_{k=1}^K  \sum_{i=i_{(k-1)}}^{i_{(k)}-1} 1_{ [T_k- \varepsilon ,  T_k)}(t_i) (t_{i+1}-t_i)(T_k-t_i)
\leqq K \varepsilon^2,
$$
So we have Equation (\ref{sum1}).\\
Next we show Equation (\ref{sum2}).
\begin{align*}
&\sum_{i=0}^{n-1} (t_{i+1}-t_i) (\sum_{k_1; T_{k_1} \geqq t_{i+1}} \phi^{(k_1)}(t_i, \varepsilon; a_1, \alpha_1, b_1, \beta_1)) 
 (\sum_{k_1; T_{k_1} \geqq t_{i+1}} \phi^{(k_2)}(t_i, \varepsilon; a_2, \alpha_2, b_2, \beta_2))  \\
\leqq &  \sum_{k=1}^K \sum_{i=i_{(k-1)}}^{i_{(k)}-1} (t_{i+1}-t_i) \sum_{j=1}^4 I^{(k)}_{i,j},
\end{align*}
where
$$
I^{(k)}_{i,1} = \sum_{k_1, k_2 = k}^K a_1 a_2 (T_{k_1} -t_i)^{-\alpha_1} (T_{k_2} -t_i)^{-\alpha_2} 1_{[0, T_{k_1} - \varepsilon)} (t_i) 1_{[0, T_{k_2} - \varepsilon)} (t_i),
$$
$$
I^{(k)}_{i,2} = \sum_{k_1, k_2 = k}^K a_1 b_2 (T_{k_1} -t_i)^{-\alpha_1} (T_{k_2} -t_i)^{\beta_2} 1_{[0, T_{k_1} - \varepsilon)} (t_i) 1_{[T_{k_2} - \varepsilon, T_{k_2})} (t_i),
$$
$$
I^{(k)}_{i,3} = \sum_{k_1, k_2 = k}^K a_2 b_1 (T_{k_2} -t_i)^{-\alpha_2} (T_{k_1} -t_i)^{\beta_1} 1_{[0, T_{k_2} - \varepsilon)} (t_i) 1_{[T_{k_1} - \varepsilon, T_{k_1})} (t_i),
$$
$$
I^{(k)}_{i,4} = \sum_{k_1, k_2 = k}^K b_1 b_2 (T_{k_1} -t_i)^{\beta_1} (T_{k_2} -t_i)^{\beta_2} 1_{[T_{k_1} - \varepsilon, T_{k_1})} (t_i)1_{[T_{k_2} - \varepsilon, T_{k_2})} (t_i).
$$
Note that (\ref{discEps}) and
\begin{align*}
1_{[0, T_{k_1} - \varepsilon)} (t_i) 1_{[T_{k} - \varepsilon, T_{k})} (t_i)=
\begin{cases}
0, k_1\leqq k\\
1_{[T_{k} - \varepsilon, T_{k})} (t_i), k_1 >k,
\end{cases}
\end{align*}
we have for $\ i \in \{i_{(k-1)}, \ldots, i_{(k)-1}\}$,
$$
I^{(k)}_{i,2} = \sum_{k_1=k+1}^K a_1 b_2 (T_{k_1} -t_i)^{-\alpha_1} (T_{k} -t_i)^{\beta_2} 1_{[T_{k} - \varepsilon, T_{k})} (t_i)
$$
$$
\leqq K(T_{k+1}-T_k)^{-\alpha_1} a_1 b_2 (T_{k} -t_i)^{\beta_2} 1_{[T_{k} - \varepsilon, T_{k})} (t_i).
$$
We have the followings similarly.
\begin{align*}
&\sum_{i=0}^{n-1} (t_{i+1}-t_i) I^{(k)}_{i,1} \leqq a_1 a_2 \sum_{k_1, k_2 = k}^K \int_{t_{i_{(k-1)}}\wedge (T_{k_1} \wedge T_{k_2}-\varepsilon)}^{t_{i_{(k)-1}} \wedge (T_{k_1} \wedge T_{k_2}-\varepsilon)} (T_{k_1}\wedge T_{k_2} -t)^{-(\alpha_1+\alpha_2)} dt,\\
&\sum_{i=0}^{n-1} (t_{i+1}-t_i) I^{(k)}_{i,2} \leqq C a_1 b_2 \varepsilon^{\beta_2 +1},\\
&\sum_{i=0}^{n-1} (t_{i+1}-t_i) I^{(k)}_{i,3} \leqq C a_2 b_1 \varepsilon^{\beta_1+1},\\
&\sum_{i=0}^{n-1} (t_{i+1}-t_i) I^{(k)}_{i,4} \leqq b_1 b_2 \varepsilon^{\beta_1+\beta_2 +1}.
\end{align*}
So we obtain (\ref{sum2}).
\qed 

\section{Proof of Theorem \ref{main2} and Theorem \ref{main1}}
\begin{thm} \label{convergence}
There exists a constant $C >0$ such that
$$E^P[|\hat{c}_1(\varepsilon_L,\Delta,L)-c_{\Delta}|] \leqq C\left( L^{-(1-\delta)/2}  \hat{e}\left(\varepsilon, (1+\delta)(\tilde{N}+1)\ell_0/4 \right) + \varepsilon^2\right), L\geqq 1.$$
\end{thm}
{\it Proof.}
$$E^P[|\hat{c}_1(\varepsilon_L,\Delta,L)-c_{\Delta}|]$$
$$ \leqq \frac{1}{L}\sum_{\ell=1}^L \sum_{i=0}^{n-1} (t_{i+1}-t_i) \sum_{m=1}^M \sum_{k:T_k\geqq t_{i+1}} |E^P[g(t_i,X(t_i,x))$$
$$((Q_{t_i,T_k,\varepsilon}^{(m,L,\omega)}F_{m,k})(\pi_k (X_{\ell}(t_i))- (P_{T_k-t_i}^{(k)}F_{m,k})(\pi_k (X_{\ell}(t_i)))]|.$$
Then by Schwartz's inequality ,
$$\frac{1}{L}\sum_{\ell=1}^L |E^P[g(t_i,X(t_i,x))
(Q_{t_i,T_k,\varepsilon}^{(m,L,\omega)}F_{m,k})(\pi_k (X_{\ell}(t_i))
- (P_{T_k-t_i}^{(k)}F_{m,k})(\pi_k (X_{\ell}(t_i))]|$$
$$\leqq \frac{1}{L}\sum_{\ell=1}^L E^P[|g(t_i,X(t_i,x))||
(Q_{t_i,T_k,\varepsilon}^{(m,L,\omega)}F_{m,k})(\pi_k (X_{\ell}(t))
- (P_{T_k-t_i}^{(k)}F_{m,k})(\pi_k (X_{\ell}(t))|1_{[0, T_k-\varepsilon)}]$$
$$+ C(T_k-t)1_{[T_k-\varepsilon, T_k)}$$
$$\leqq \left( \frac{1}{L}\sum_{\ell=1}^L E^P[|
 (Q_{t_i,T_k,\varepsilon}^{(m,L,\omega)}F_{m,k})(\pi_k (X_{\ell}(t))
- (P_{T_k-t_i}^{(k)}F_{m,k})(\pi_k (X_{\ell}(t))|^2] \right)^{1/2} 1_{[0, T_k-\varepsilon)}$$
$$\times E[g(t_i,X(t_i,x))^2]^{1/2} + C(T_k-t)1_{[T_k-\varepsilon, T_k)}$$
By Proposition \ref{key}
$$E^P[|\hat{c}_1(\varepsilon_L,\Delta,L)-c_{\Delta}|]$$
$$\leqq C \sum_{i=0}^{n-1} (t_{i+1}-t_i) \sum_{k:T_k\geqq t_{i+1}}
 \phi^{(k)}(t_i, \varepsilon;L^{-(1-\delta)/2}, (1+\delta)(\tilde{N}+1)\ell_0/4, 1, 1).$$
By Proposition\ref{dd},
we have the assertion. 
\qed

\begin{lem} \label{improve}
Let $a, b \in {\bf R} $ and $c, \theta >0$. Then we have
$$c|a| |1_{\{b \geqq 0\}} -1_{\{a \geqq 0\}}| \leqq c|b-a| 1_{\{ |b-a| \geqq \theta \}} +c\theta 1_{\{ |a| < \theta \}} $$
\end{lem}
{\it Proof.}
If $|a| > |a-b|$, then 
$$ 1_{\{b \geqq 0\}} -1_{\{a \geqq 0\}} = 0.$$
So we see that
$$ |a|( |1_{\{b \geqq 0\}} -1_{\{a \geqq 0\}}|$$
$$\leqq |a| 1_{\{ |a| \leqq |a-b|\}}$$
$$\leqq |a-b|1_{\{|a-b| \geqq \theta \}} +|a|1_{\{ |a| < \theta \}}.$$
\qed

\begin{thm}
Let $\delta \in (0,1), p>1.$ 
Suppose that there is $ \gamma \in (0,1]$  and $C_{\gamma} >0$ such that
  $$\sup_{\Delta} \sum_{i=0}^{n-1} (t_{i+1}-t_i)){\mu} (| \sum_{m=1}^M\sum_{k:T_k\geqq t_{i+1}} (P_{T_k-t_i}^{(m)}F_{m,k})(\pi_m X(t_i,x^*)))| \leqq \theta )$$
$$ \leqq C_{\gamma}\theta^{\gamma}, \theta\in(0,1).$$
Then there exists a constant $C>0$ , $\tilde{\Omega}(L, \varepsilon) \in \mathcal{F}$, such that
$$P( \tilde{\Omega}(L,\varepsilon)^c ) \leqq C \left( \left(
L^{-(1-\delta)^2/2} \left( L^{-(1-\delta)^2/2} \hat{e}\left(\varepsilon,(1-\delta^2)(\tilde{N}+1)\ell_0/2 \right)+\varepsilon^{3(1-\delta)/2} \right)
  \right)^{\frac{ \delta(1+\gamma)}{2+\gamma}} \right),$$
and
$$
1_{\tilde{\Omega}(L,\varepsilon)}|\hat{c}_2(\varepsilon_L,\Delta,L)-c_{\Delta}|$$
$$\leqq  C\left(
L^{-(1-\delta)^2/2} \left( L^{-(1-\delta)^2/2} \hat{e}\left(\varepsilon,(1-\delta^2)(\tilde{N}+1)\ell_0/2 \right)+\varepsilon^{3(1-\delta)/2} \right)
  \right)^{\frac{(1-\delta)(1+\gamma)}{2+\gamma}}, $$
  for $ L \geqq 1.$
\end{thm}
{\it Proof.}
In this proof, we denote $X(t,x^*)$ by $X(t)$ for simplicity.
Let $$\tilde{B}_{\varepsilon} = \bigcap_{m=1}^M \bigcap_{k=1}^K B^{(m,k)}(T_k-\varepsilon, \delta,L).$$
Let $F_{P,i } : {\bf R}^N \to {\bf R}$ be given by
$$F_{P,i }(x) = \sum_{\substack{m=1,\ldots,M \\ k:T_k\geqq t_{i+1}}} (P_{T_k-t_i}^{(m)} F_{m,k}) (\pi_m x),$$
and let $F_{Q,i } : {\bf R}^N \to {\bf R}$ be given by
$$F_{Q,i }(x)=\sum_{\substack{m=1,\ldots,M \\ k:T_k\geqq t_{i+1}} } (Q_{t_i,T_k, \varepsilon}^{(m)}F_{m,k})(\pi_m x).$$
Then
$$1_{\tilde{B}_{\varepsilon} }|\hat{c}_2(\varepsilon_L,\Delta,L)-c_{\Delta}|$$
 $$
  \leqq 1_{\tilde{B}_{\varepsilon} }\sum_{i=0}^{n-1}(t_{i+1}-t_i)
|E ^{\mu}[|g(t_i, X(t_i))|\sum_{m=1}^M \sum_{k:T_k\geqq t_{i+1}}F_{m,k} (\pi_k (X(T_k)))
(1_{\{F_{Q,i}(X(t_i) \geqq 0\}} -1_{\{F_{P,i}(X(t_i) \geqq 0\}})]|,$$
  \begin{align*} 
   \leqq &1_{\tilde{B}_{\varepsilon} } \sum_{i=0}^{n-1}(t_{i+1}-t_i)
E ^{\mu}[|g(t_i, X(t_i))||F_{P,i}(X(t_i)| 
|1_{\{F_{Q,i}(X(t_i) \geqq 0\}} -1_{\{F_{P,i}(X(t_i) \geqq 0\}}|],
 \end{align*}
since $|g(t_i, X(t_i))|(1_{\{F_{Q,i}(X(t_i) \geqq 0\}} -1_{\{F_{P,i}(X(t_i) \geqq 0\}})$ is $\mathcal{F}_t$ measurable.
 \\
Applying Lemma \ref{improve} to $a = F_{P,i}(X(t_i), b = F_{Q,i}(X(t_i),$ and $c=|g(t_i, X(t_i))|,$ we have
$$1_{\tilde{B}_{\varepsilon} }|\hat{c}_2(\varepsilon_L,\Delta,L)-c_{\Delta}|\leqq I_1 +I_2, $$
where
\begin{align*}
&I_1=1_{\tilde{B}_{\varepsilon} } \sum_{i=0}^{n-1}(t_{i+1}-t_i)
E ^{\mu}[|g(t_i, X(t_i))||F_{Q,i}(X(t_i) - F_{P,i}(X(t_i)| 1_{\{|F_{Q,i}(X(t_i) - F_{P,i}(X(t_i)| > \theta \}}],\\
&I_2=1_{\tilde{B}_{\varepsilon} } \sum_{i=0}^{n-1}(t_{i+1}-t_i) \theta E^{\mu} [|g(t_i, X(t_i))||1_{\{F_{P,i}(X(t_i)| \leqq \theta \}}].
\end{align*}
By H\"{o}lder's inequality,
$$I_2 \leqq \theta E^{\mu}[|g(t_i, X(t_i))|^{1/\delta}]^{\delta} \mu (|F_{P,i}(X(t_i)| \leqq \theta)^{1-\delta} \leqq  C \theta \mu (|F_{P,i}(X(t_i)| \leqq \theta).$$
From the assumption, we have
$$I_2 \leqq C \theta^{(\gamma+1)(1-\delta)}.$$
Next we will estimate $I_1$.
$$
|F_{Q,i}(X(t_i)) - F_{P,i}(X(t_i))| \leqq \sum_{m=1}^M\sum_{k:T_k\geqq t_{i+1}}  (d_{1,\varepsilon,L}^{(m,k)}(t_i,X(t_i))
+ d_{2,\varepsilon,L}^{(m,k)}(t_i,X(t_i)) +d_{3,\varepsilon,L}^{(m,k)}(t_i,X(t_i)))
$$
\begin{align*}
&I_1\leqq 1_{\tilde{B}_{\varepsilon} } \sum_{i=0}^{n-1}(t_{i+1}-t_i)  E ^{\mu}[ |g(t_i,X(t_i))| \\
&\sum_{m=1}^M \sum_{k:T_k\geqq t_{i+1} }\left(d_{1,\varepsilon,L}^{(m,k)}(t_i, \pi_m (X(t_i))) +d_{2,\varepsilon,L}^{(m, k)}(t_i, \pi_m (X(t_i) + d_{3,\varepsilon,L}^{(m, k)}(t_i, \pi_m (X(t_i) \right)  \\
&\times 1_{\{\sum_{m=1}^M\sum_{k:T_k\geqq t_{i+1} }  (d_{1,\varepsilon,L}^{(m, k)}(t_i, \pi_m (X(t_i))) +d_{2,\varepsilon,L}^{(m, k)}(t_i, \pi_m (X(t_i)) + d_{3,\varepsilon,L}^{(m, k)})(t_i, \pi_m (X(t_i))) > \theta\} }]  .
\end{align*}
$I_1$ is dominated by
$$I_1\leqq I_{1,1}+I_{1,2}+I_{1,3}+I_{1,4},$$
where
$$
I_{1,1}=\sum_{i=0}^{n-1}(t_{i+1}-t_i)\sum_{m=1}^M \sum_{k:T_k\geqq t_{i+1} } E ^{\mu}[|g(t_i,X(t_i))|d_{1,\varepsilon,L}^{(m,k)}(t_i, \pi_m (X(t_i))) ],
$$
$$
I_{1,2}= \sum_{i=0}^{n-1}(t_{i+1}-t_i)   E ^{\mu}[\sum_{m=1}^M \sum_{k:T_k\geqq t_{i+1} } |g(t_i,X(t_i))|d_{2,\varepsilon,L}^{(m,k)}(t_i, \pi_m (X(t_i)))$$
$$\times 1_{\{\sum_{m=1}^M \sum_{k:T_k\geqq t_{i+1} } (d_{1,\varepsilon,L}^{(m,k)}(t_i, \pi_k (X(t_i) + d_{3,\varepsilon,L}^{(m,k)}(t_i, \pi_m (X(t_i))  > \theta/2\}}],
$$
$$
I_{1,3}= \sum_{i=0}^{n-1}(t_{i+1}-t_i)  E ^{\mu}[\sum_{m=1}^M \sum_{k:T_k\geqq t_{i+1} } |g(t_i,X(t_i))|d_{2,\varepsilon,L}^{(m,k)}(t_i, \pi_m (X(t_i) ))$$
$$\times 1_{\{\sum_{m=1}^M \sum_{k:T_k\geqq t_{i+1} } d_{2,\varepsilon,L}^{(m,k)}(t_i, \pi_m (X(t_i)  > \theta/2\}}].
$$
$$I_{1,4}= \sum_{i=0}^{n-1}(t_{i+1}-t_i)  E ^{\mu}[\sum_{m=1}^M \sum_{k:T_k\geqq t_{i+1} } |g(t_i,X(t_i))|d_{3,\varepsilon,L}^{(m,k)}(t_i, \pi_m (X(t_i) ))]$$
From Proposition \ref{d12},
$$
E^P[I_{1,1}]=\sum_{i=0}^{n-1} (t_{i+1}-t_i) \sum_{\substack{m=1,\ldots,M,\\ k; T_k \geqq t_{i+1}} } \int_{E} E^P[d_{1,\varepsilon,L}^{(m,k)}(t,\tilde{x}_m) ]
|g(t_i,x)| p(t_i, x^*, x)dx 
$$
$$
\leqq C\sum_{i=0}^{n-1} (t_{i+1}-t_i) \sum_{k; T_k \geqq t_{i+1}}  \phi^{(k)}(t_i,\varepsilon; L^{-(1-\delta)^3}, 0,0, 0)$$
$$
\leqq C\hat{e}\left(\varepsilon,L^{-(1-\delta)^3}\right).
$$
\\
Next, we will estimate $I_{1,2}$.  By H\"{o}lder's inequality
$$I_{1,2}\leqq 
 \sum_{i=0}^{n-1}(t_{i+1}-t_i)   E ^{\mu}[(\sum_{m=1}^M \sum_{k:T_k\geqq t_{i+1} } d_{2,\varepsilon,L}^{(m, k)}(t_i, \pi_m (X(t_i)))^2)]^{1/2} 
 $$
 $$
\times E^{\mu}[|g(t_i,X(t_i))|^{2/\delta}]^{\delta/2}E^{\mu}[1_{\{\sum_{m=1}^M \sum_{k:T_k\geqq t_{i+1} } (d_{1,\varepsilon,L}^{(m,k)}(t_i, \pi_m(X(t_i) + d_{3,\varepsilon,L}^{(m,k)}(t_i, \pi_m (X(t_i) ) > \theta/2\}}]^{(1-\delta)/2},$$
$$
 \leqq C\sum_{i=0}^{n-1}(t_{i+1}-t_i)  E ^{\mu}[\sum_{m=1}^M \sum_{k:T_k\geqq t_{i+1} } d_{2,\varepsilon,L}^{(m,k)}(t_i, \pi_m (X(t_i)))^2]^{1/2} 
 $$
 $$
\times (\frac{2}{\theta} E^{\mu}[\sum_{m=1}^M \sum_{k:T_k\geqq t_{i+1} } (d_{1,\varepsilon,L}^{(m, k)}(t_i, \pi_m (X(t_i)+ d_{3,\varepsilon,L}^{(m, k)}(t_i, \pi_m (X(t_i)) ])^{(1-\delta)/2}.$$
So we have
\begin{align*}
&E^P[ I_{1,2}] \leqq  C\sqrt{\frac{2}{\theta}} \sum_{i=0}^{n-1}(t_{i+1}-t_i)  
E^P[\sum_{\substack{m=1,\ldots,M,\\k:T_k\geqq t_{i+1}} }  \int_{E_m} d_{2,\varepsilon,L}^{(m,k)}(t_i, x)^2 p^{(m)}(t_i, \tilde{x}_m^*, x)dx ]^{1/2}\\
&\times E^P[ \sum_{\substack{m=1,\ldots,M,\\k:T_k\geqq t_{i+1}} } \int_{E_m} (d_{1,\varepsilon,L}^{(m,k)}(t_i, x) + d_{3,\varepsilon,L}^{(m,k)}(t_i, x)) p^{(m)}(t_i, \tilde{x}_m^*, x)dx ]^{(1-\delta)/2}\\
\leqq &C \theta^{-(1-\delta)/2} \sum_{i=0}^{n-1}(t_{i+1}-t_i)  \left( \sum_{k:T_k\geqq t_{i+1} } \phi^{(k)}(t_i, \varepsilon; L^{-(1-\delta)^2/2}, (1-\delta^2)(\tilde{N}+1)\ell_0/4, 0, 0)\right) \\
&\times \left(\sum_{k:T_k\geqq t_{i+1} } \phi^{(k)}(t_i, \varepsilon; L^{-(1-\delta)^4/2}, 0, 1, (1-\delta)/2)\right).
\end{align*}
By Proposition \ref{dd},
$$
E^P[ I_{1,2}]\leqq C\theta^{-1/2}  \left( L^{-(1-\delta)^4} \hat{e}\left(\varepsilon, (1-\delta^2)(\tilde{N}+1)\ell_0/4\right)
+ L^{-(1-\delta)^2/2}\varepsilon^{3(1-\delta)/2} \right).
$$
Similarly, we have 
\begin{align*}
&E^P[ I_{1,3}] \leqq  C \theta^{-(1-\delta)/2} \sum_{i=0}^{n-1}(t_{i+1}-t_i)  
E[\sum_{\substack{m=1,\ldots,M,\\k:T_k\geqq t_{i+1}} }  \int_{E_m} d_{2,\varepsilon,L}^{(k)}(t_i, x)^2 p^{(k)}(t_i, \tilde{x}_k^*, x)dx ]^{(1-\delta)}\\
\leqq &C \theta^{-(1-\delta)/2} \sum_{i=0}^{n-1}(t_{i+1}-t_i)  (\sum_{k:T_k\geqq t_{i+1} } \phi^{(k)}(t_i, \varepsilon; L^{-(1-\delta)^2/2}, (1-\delta^2)(\tilde{N}+1)\ell_0/4, 0,0) )^2\\
  \leqq &C\theta^{-1} L^{-(1-\delta)^2} \hat{e}\left(\varepsilon,(1-\delta^2)(\tilde{N}+1)\ell_0/2 \right).
\end{align*}
It follows easily that
$$
E^P[ I_{1,4}]\leqq C\varepsilon^2.
$$
Notice that
$$
\theta^{-(1-\delta)/2}  \hat{e}(\varepsilon, (1-\delta^2)(\tilde{N}+1)\ell_0/4) \leqq
\theta^{-1}  \hat{e}\left(\varepsilon,(1-\delta^2)(\tilde{N}+1)\ell_0/2 \right),
$$
we have
$$
E^P[I] \leqq C\left(\theta^{\gamma+1}+\theta^{-1}  
L^{-(1-\delta)^2/2} \left( L^{-(1-\delta)^2/2} \hat{e}\left(\varepsilon,(1-\delta^2)(\tilde{N}+1)\ell_0/2 \right)+\varepsilon^{3(1-\delta)/2} \right) \right).
$$

In particular if we take $\theta=\theta_L$ as 
$$ 
\theta_L=O\left(
L^{-(1-\delta)^2/2} \left( L^{-(1-\delta)^2/2} \hat{e}\left(\varepsilon,(1-\delta^2)(\tilde{N}+1)\ell_0/2 \right)+\varepsilon^{3(1-\delta)/2} \right)
 \right)^{\frac{1}{2+\gamma}} ,$$
then we have
$$E^P[I] \leqq 
 C\left(
L^{-(1-\delta)^2/2} \left( L^{-(1-\delta)^2/2} \hat{e}\left(\varepsilon,(1-\delta^2)(\tilde{N}+1)\ell_0/2 \right)+\varepsilon^{3(1-\delta)/2} \right)
 \right)^{(1+\gamma)/(2+\gamma)} .$$

Let $\tilde{\Omega}(L,\varepsilon)$ be
\begin{align*}
&\tilde{\Omega}(L,\varepsilon) =\tilde{B}_{\varepsilon}\cap  \{ \omega \in \Omega; \\
&I \leqq C\left(\left(
L^{-(1-\delta)^2/2} \left( L^{-(1-\delta)^2/2} \hat{e}\left(\varepsilon,(1-\delta)^2(\tilde{N}+1)\ell_0/2 \right)+\varepsilon^{3(1-\delta)/2} \right)
 \right)
  \right)^{\frac{(1-\delta)(1+\gamma)}{2+\gamma}}  \}.
\end{align*}
From Proposition \ref{exSet}, we have
$$P( \tilde{\Omega}(L,\varepsilon)^c ) \leqq C \left( \left(
L^{-(1-\delta)^2/2} \left( L^{-(1-\delta)^2/2} \hat{e}\left(\varepsilon,(1-\delta^2)(\tilde{N}+1)\ell_0/2 \right)+\varepsilon^{3(1-\delta)/2} \right)
  \right)^{\frac{ \delta(1+\gamma)}{2+\gamma}} \right.$$
$$\left.   + (\varepsilon^{-5\ell_0}L^{-p\delta^2/2+1/p})^p\right),$$
and
$$
1_{\tilde{\Omega}(L,\varepsilon)}|\hat{c}_2(\varepsilon_L,\Delta,L)-c_{\Delta}|$$
$$\leqq  C\left(
L^{-(1-\delta)^2/2} \left( L^{-(1-\delta)^2/2} \hat{e}\left(\varepsilon,(1-\delta^2)(\tilde{N}+1)\ell_0/2 \right)+\varepsilon^{3(1-\delta)/2} \right)
  \right)^{\frac{(1-\delta)(1+\gamma)}{2+\gamma}} .$$
 \qed
 
\begin{cor}

Theorem \ref{main2} and
Theorem \ref{main1} follow from the Theorem \ref{convergence} and Theorem \ref{improve}.

\end{cor}
 
\section{Numerical Example}
\begin{figure}[htbp]
\begin{center}
\includegraphics[width=12cm,bb=0 0 500 400]{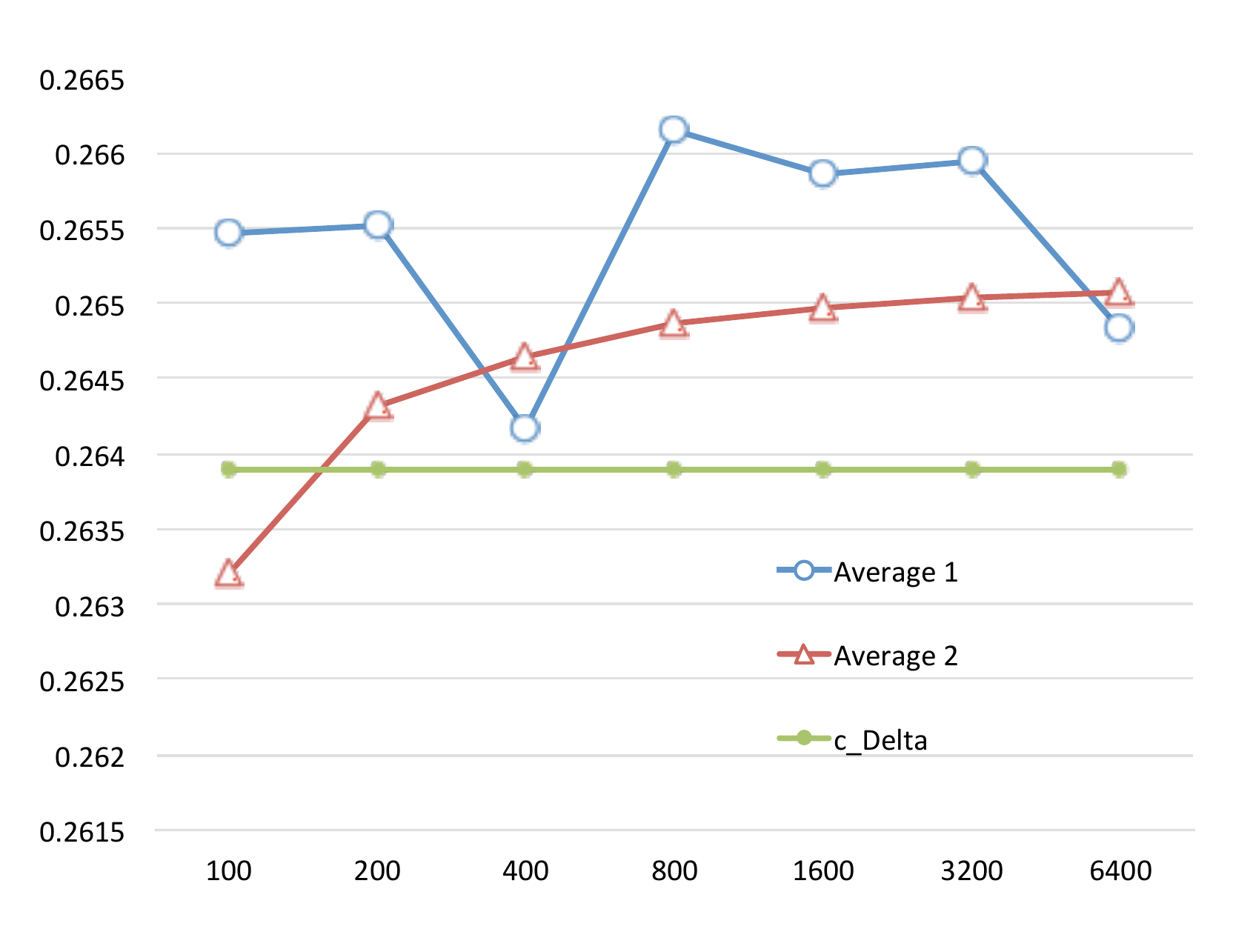}
\caption{Average} \label{fig1}
\end{center}
\end{figure}

\begin{figure}[htbp]
\begin{center}
\includegraphics[width=12cm,bb=0 0 500 400]{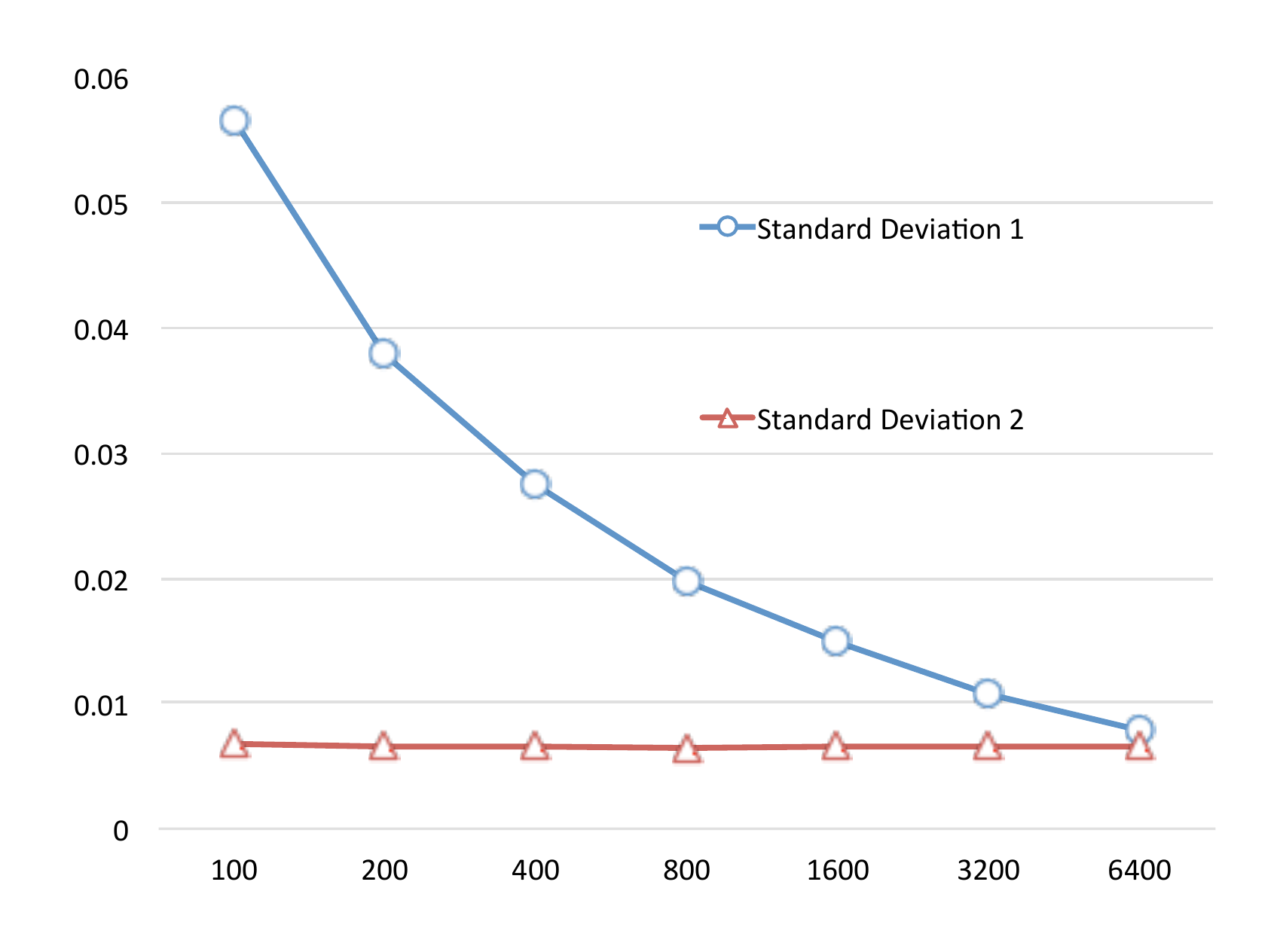}
\caption{Standard Deviation} \label{fig2}
\end{center}
\end{figure}

Let $\{B(t); t\geqq 0\} $ be 1 dimensional Brownian motion.
Let $t_i = i/n, i=0,\ldots, n.$ Let $c$ be
$$c=E [\int_0^1( E[B(1)|\mathcal{F}_t]\vee0 )\ dt ]=\frac{2}{3 \sqrt{2\pi}}.$$
Let $c_{\Delta}$ be the discretization of $c$, such that
$$
c_{\Delta}=\sum_{i=0}^{n-1}(t_{i+1}-t_i)E[B(t_i)\vee0].
$$
We approximate $c$  as Remark \ref{actual}, where $F(x) =x.$ Let
${\bf X}_1 = \{X_{\ell }(t_i); i=0,1,\ldots,n \}_{\ell=1}^L$
be i.i.d sample paths of $\{B(t_i); i=0,1,\ldots,n\} $.
We compute
 $Q_{t,T,\varepsilon}^{(L,\omega)}$ and $\hat{c}_1$ by using of paths ${\bf X}_1$.
$$\hat{c}_1=\frac{1}{L} \sum_{\ell=1}^L \sum_{i=0}^{n-1} (Q_{t_i,T}^{(L,\omega)}F)((X_{\ell}(t_i)))\vee0)(t_{i+1}-t_i). $$
Let ${\bf X}_2 = \{X'_{\ell }(t_i); i=0,1,\ldots,n \}_{\ell=1}^L$
be another i.i.d sample paths of $\{B(t_i); i=0,1,\ldots,n\} $.
We compute $\hat{c}_2$ by 
 $$\hat{c}_2= \frac{1}{L_0} \sum_{\ell_0=1}^{L_0} \sum_{i=0}^{n-1} \{ F(X_{\ell_0}^{'}(T))\}
 1_{\{(Q_{t,T}^{(L,\omega)}F)(X^{'}_{\ell_0}(t_i)) \geqq 0\}}(t_{i+1}-t_i).$$ 
We have $c\approx 0.2659615203.$ When we take $n = 100,$ we have $c_{\Delta}\approx 0.2638855365.$ 
We also take $L_0=10000$ and $L=100, 200, 400, 800, 1600, 3200, 6400.$ 
We replicate 100 estimators of $\hat{c}_i, i=1,2$ for each $L$. Let "Average $i$" denote the average and "Standard Dviation $i$" denote
 the unbiased standard deviation of these 100 estimators of $\hat{c}_i, i=1,2.$ We show the numerical result in Table \ref{nume}, 
we show graph of "Average $i, i = 1,2$" and $c_{\Delta}$ in Figure \ref{fig1} and graph of "Standard Deviation $i, i=1,2$" in Figure \ref{fig2}.
\noindent 
We see in Figure \ref{fig1} that both $\hat{c}_1$ and $\hat{c}_2$ are close to $c_{\Delta}$,
but we see in Figure \ref{fig2} that $\hat{c}_2$ is more stable than $\hat{c}_1.$
\\
\begin{table}[htb]
\begin{center} 
  \begin{tabular}{ccccc}
    $L$ & Average 1& Average 2 &  Standard Deviation 1 & Standard Deviation 2\\ \hline 
	100	&	0.2654599783 	&	0.2632060491 	&	0.0334099140 	&	0.0066656532 	\\
	200	&	0.2655180439 	&	0.2643180950 	&	0.0244301557 	&	0.0064812199 	\\
	400	&	0.2641632386 	&	0.2646417384 	&	0.0168528548 	&	0.0064796412 	\\
	800	&	0.2661557058 	&	0.2648673390 	&	0.0114840865 	&	0.0064426983 	\\
	1600	&	0.2658622710 	&	0.2649728715 	&	0.0092734158 	&	0.0064583976 	\\
	3200	&	0.2659440890 	&	0.2650330723 	&	0.0071782226 	&	0.0064616852 	\\
	6400	&	0.2648318867 	&	0.2650702303 	&	0.0050607508 	&	0.0064608708 	\\

  \end{tabular}
  \caption{Average and Standard Deviation} \label{nume}
\end{center}
\end{table}

\newpage
\begin{thebibliography}{99}
\bibitem{AH}Avramidis, A.N. , and P. Hyden, 
\newblock Efficiency improvements for pricing American options 
with a stochastic mesh, 
\newblock in Proceedings of the 1999 Winter Simulation Conference, pp. 344-350
\bibitem{BG1} Broadie, M., and P. Glasserman,
\newblock A stochastic mesh method for pricing high-dimensional American options
\newblock J. Computational Finance, 7 (4) (2004), 35-72.
\bibitem{DU} Duffie, D., M. Huang, 1996,
\newblock  Swap Rates and Credit Quality, 
\newblock Journal of Finance, Vol. 51, No. 3, 921
\bibitem{FT} Fujii, F., A. Takahashi
\newblock Derivative Pricing under Asymmetric and Imperfect Collateralization, and CVA
\newblock CARF Working Paper F-240, December 2010
\bibitem{G}  Glasserman, P.,
\newblock "Monte Carlo Methods in Financial Engineering"
\newblock Springer,  2004, Berlin. 
\bibitem{Gre}Gregory, J.,
\newblock Counterparty credit risk and credit value adjustment: A continuing challenge for global financial markets, Second Editio,
\newblock John Wiley \& Sons, 2012.
\bibitem{K}Kusuoka, S.,
\newblock Malliavin Calculus Revisited,
\newblock  J. Math. Sci. Univ. Tokyo 10(2003), 261-277.
\bibitem{KS2}Kusuoka, S., and D.W.Stroock, 
\newblock Applications of Malliavin Calculus II,
\newblock J. Fac. Sci. Univ. Tokyo Sect. IA Math. 32(1985),1-76.
\bibitem{KM}Kusuoka, S. , and Y. Morimoto, 
\newblock Stochastic mesh methods for H{\" o}rmander type diffusion processes, 
\newblock Preprint.
\bibitem{H} Laborde're,H.
\newblock Cutting CVA's Complexity,
\newblock Risk Magazine, July, 2012, pp. 67-73.
\bibitem{LH} Liu, G., and L.J. Hong,
\newblock Revisit of stochastic mesh method for pricing American options,
\newblock Operations Research Letters 37(2009), 411-414.
\bibitem{Sh}Shigekawa, I.,
\newblock "Stochastic Analysis",
\newblock Translation of Mathematical Monographs vol.224,
AMS 2000.
\end{thebibliography}
\end{document}